\documentclass{article}
\usepackage[utf8]{inputenc}
\usepackage{enumitem}
\usepackage{siunitx}

\usepackage{pgfplots}
\pgfplotsset{compat=1.8}
\usepackage{pgfplotstable}

\title{The Smallest Hard Trees}

\author{Manuel Bodirsky\footnote{Institut f\"ur Algebra, TU Dresden. https://orcid.org/0000-0001-8228-3611.}, Jakub Bul\'in\footnote{Department of Theoretical Computer Science and Mathematical Logic, Faculty of Mathematics and Physics, Charles University, Prague. https://orcid.org/0000-0001-5235-8715.}, Florian Starke\footnote{Institut f\"ur Algebra, TU Dresden. https://orcid.org/0000-0003-2360-9364.}, Michael Wernthaler\footnote{Institut f\"ur Algebra, TU Dresden. https://orcid.org/0000-0003-4978-0928.}}

\date{May 2022}

\usepackage{amsmath,amssymb,a4wide,url,graphicx}
\usepackage[english]{babel}
\usepackage{amsxtra,amssymb,amsthm,hyperref}
\usepackage[linesnumberedhidden]{algorithm2e}
\usepackage{makecell,booktabs}
\usepackage{tikz}
\tikzstyle{bullet}=[circle,fill,inner sep=0pt,minimum size=3pt]
\tikzstyle{bullethole}=[circle,draw,inner sep=0pt,minimum size=3pt]
\usetikzlibrary{arrows}

\newcommand{\ignore}[1]{}
\DeclareMathOperator{\Csp}{CSP}
\DeclareMathOperator{\dist}{dist}
\DeclareMathOperator{\id}{id}
\DeclareMathOperator{\AC}{AC}
\DeclareMathOperator{\End}{End}

\DeclareMathOperator{\levelfunction}{lvl}

\newcommand{\structure}[1]{\mathbb{#1}}
\newcommand{\digraph}[1]{\structure{#1}}
\newcommand{\edges}[1]{E({#1})}
\newcommand{\level}[1]{\levelfunction({#1})}

\newcommand{\G}{\digraph G}
\newcommand{\T}{\digraph T}
\newcommand{\HH}{\digraph H}
\newcommand{\indicator}[1]{{#1}^\mathrm{Ind}}
\newcommand{\edgesH}{\edges{\HH}}

\newcommand{\TLinP}{\ensuremath{{\mathrm{3Lin}}_p}}

\DeclareMathOperator{\HornSAT}{Horn-3SAT}
\DeclareMathOperator{\stcon}{st-Con}

\newtheorem{theorem}{Theorem}
\newtheorem{question}{Question}
\newtheorem{conjecture}{Conjecture}

\newtheorem{corollary}[theorem]{Corollary}
\newtheorem{lemma}[theorem]{Lemma}
\newtheorem{remark}[theorem]{Remark}
\newtheorem{definition}[theorem]{Definition}

\DeclareMathOperator{\KK}{KK}
\DeclareMathOperator{\HMcK}{HMcK}
\DeclareMathOperator{\HM}{HM}
\DeclareMathOperator{\NN}{NN}
\DeclareMathOperator{\TS}{TS}
\DeclareMathOperator{\J}{J}
\DeclareMathOperator{\NU}{NU}

\newcommand{\blue}[1]{
\textcolor{blue}{#1}}

\newcommand\keywordslabel{Keywords:}
\newcommand\keywords[1]{%
  \begin{list}{}{%
    \setlength{\topsep}{2ex}%
    \settowidth{\leftmargin}{\bfseries\keywordslabel~}%
    \setlength{\labelsep}{0pt}%
    \setlength{\labelwidth}{\leftmargin}%
    \setlength{\itemindent}{0pt}%
  }
  \raggedright\item[\bfseries\keywordslabel~]#1
  \end{list}
}

\begin{document}

\maketitle

\begin{abstract}
\noindent We find an orientation of a tree with 20 vertices such that the corresponding fixed-template constraint satisfaction problem (CSP) is NP-complete, and prove that for every orientation of a tree with fewer vertices the corresponding CSP can be solved in polynomial time. 
We also compute the smallest tree that is NL-hard (assuming $\operatorname{L}\neq\operatorname{NL}$), the smallest tree that cannot be solved by arc consistency, and the smallest tree that cannot be solved by Datalog. 
Our experimental results also support a conjecture of Bul\'in concerning a question of Hell, Ne\v{s}et\v{r}il and Zhu, namely that `easy trees lack the ability to count'.
Most proofs are computer-based and make use of the most recent universal-algebraic theory about the complexity of finite-domain CSPs. However, further ideas are required  because of the huge number of orientations of trees. In particular, we use the well-known fact that it suffices to study orientations of trees that are cores and show how to efficiently decide whether a given orientation of a tree is a core using the arc-consistency procedure. Moreover, we present a method to generate orientations of trees that are cores that works well in practice. In this way we found interesting examples for the open research problem to classify finite-domain CSPs in NL. 
\\

\end{abstract}

\section{Introduction}
For a fixed directed graph (short: \emph{digraph}) $\HH$, the \emph{constraint satisfaction problem for $\HH$}, denoted by $\Csp(\HH)$, is the computational problem of deciding whether a given finite digraph $\G$ admits a homomorphism to $\HH$. This problem is also known as the \emph{$\HH$-coloring problem}.
If $\HH$ is finite and symmetric, then $\Csp(\HH)$ can be solved in polynomial time if $\HH$ is bipartite or contains a loop, and is NP-complete otherwise~\cite{HellNesetril}.  The situation if $\HH$ is a finite but not necessarily symmetric digraph is much more complicated. 
The Feder-Vardi dichotomy conjecture states that $\Csp(\HH)$ is in P or NP-complete~\cite{FederVardi}. 
In fact, the conjecture was phrased not only for digraphs but for the corresponding computational problem for general finite relational structures. 
However, it is known that every CSP for a finite relational structure is polynomial-time (and even logspace) 
equivalent to the CSP for a finite digraph~\cite{FederVardi,BulinDelicJacksonNiven}. 
The Feder-Vardi conjecture was proved in 2017 independently by Bulatov~\cite{BulatovFVConjecture} and by Zhuk~\cite{ZhukFVConjecture,Zhuk20}. Prior to their break-through result, the conjecture was open even if $\HH$ is an orientation of a finite tree.

The description of the polynomially solvable cases in the proofs of Bulatov and of Zhuk is based on the so-called algebraic approach 
and phrased using  \emph{polymorphisms} of $\HH$, i.e., edge-preserving multivariate operations on the vertex set (`higher-dimensional symmetries') \cite{BartoKozikWillardPolymorphisms}. The algebraic condition for polynomial-time tractability in the proofs of Bulatov and of Zhuk has numerous equivalent characterizations, e.g. \cite{Siggers,wnuf,Cyclic}. Siggers was the first to show the (at the time somewhat surprising) fact that the condition can be characterized by the existence of a single, 6-ary polymorphism satisfying certain identities \cite{Siggers} --- which can readily be tested at least for very small digraphs. This was later improved by Kearnes, Markovi\'c, and McKenzie \cite{KearnesMarkovicMcKenzie} to the existence of a single 4-ary operation, commonly referred to as a \emph{Siggers polymorphism}, or a pair of 3-ary operations which we will call \emph{Kearnes-Markovi\'{c}-McKenzie polymorphisms} (for the definition, see Section~\ref{sect:Siggers}). The latter is computationally the most feasible (the search space is the smallest) and thus the most suitable for our purposes. The question whether a given finite digraph $\HH$ satisfies any of the equivalent characterizations of the algebraic tractability condition is decidable, but NP-hard~\cite{MetaChenLarose}.

\subsection{Computational Complexity}
Several other important  conjectures about the computational complexity of the constraint satisfaction problem for a fixed finite structure $\HH$ with finite relational signature  remain open: most notably the question for which finite structures $\HH$ the problem $\Csp(\HH)$ is in the complexity class NL (non-deterministic logspace), and for which finite structures $\HH$ it is in the complexity class L (deterministic logspace). 
As in the case of P versus NP-hard, it appears that these questions are closely linked to central dividing lines in universal algebra, as illustrated by the following conjectures.

\begin{conjecture}[Larose and Tesson~\cite{LaroseTesson}]\label{conj:NL} If the polymorphisms of a finite structure $\HH$ with finite relational signature contain a Kearnes-Kiss chain (defined in Section~\ref{sect:PorModpL-hardness}), then $\Csp(\HH)$ is in NL. \end{conjecture}

It is known that if $\HH$ does not satisfy the condition from Conjecture~\ref{conj:NL}, then $\HH$ is hard for complexity classes that are not believed to be in NL (more details can be found in Section~\ref{sect:PorModpL-hardness}).
Conjecture~\ref{conj:NL}
is wide open and we believe it to be one of the most difficult research problems in the theory of finite-domain constraint satisfaction that remains open.

\begin{conjecture}
[Egri, Larose, and Tesson~\cite{EgriLaroseTessonLogspace}]\label{conj:L}
If the polymorphisms of a finite structure $\HH$ with finite relational signature contain a Noname chain  
 (defined in  Section~\ref{sect:NLorModpL-hardness})
 then $\Csp(\HH)$ is in L.
 \end{conjecture} 
 
Also here it is known that if $\HH$ does not satisfy the condition from Conjecture~\ref{conj:L}, then $\HH$ is hard for complexity classes that are not believed to be in L (more details can be found in Section~\ref{sect:NLorModpL-hardness}).
 
We mention that both conjectures can equivalently be phrased by the inability to  \emph{primitively positively construct in $\HH$} certain finite structures that are known to be L-hard, Mod$_p$L-hard, or NL-hard (see Section~\ref{sect:specific-polymorphisms}).
Kazda proved a conditional result that states that resolving the first conjecture would also provide a solution to the second~\cite{Kazda-n-permute}. 

Again, 
these conjectures are already open if $\HH$ is a finite digraph, or even if $\HH$ is an orientation of a finite tree. It is also known that answering the question of containment in NL for finite digraphs would also answer the question for general finite structures~\cite{BulinDelicJacksonNiven}. For orientations of finite trees, however, the question might be easier to resolve.
For brevity, an orientation of a finite tree is simply called a \emph{tree} in this paper and we adopt the following terminology: a digraph $\HH$ is \emph{NP-hard} if $\Csp(\HH)$ is NP-hard, and \emph{tractable} if $\Csp(\HH)$ is in P. Similarly, we say that $\HH$ is
\emph{P-hard}, \emph{NL-hard}, \emph{in NL}, \emph{in L}, \emph{NP-complete},  \emph{NL-complete}, etc.\ if $\Csp(\HH)$ has that property.

Unfortunately, there is no graph theoretic characterization of which trees are NP-hard. The first NP-hard tree $\T$ was found by Gutjahr, Welzl, and Woeginger and had 287 vertices~\cite{GutjahrWW92}. This was later improved by Gutjahr to a smaller NP-hard tree with 81 vertices~\cite{Gutjahr}, and
then to an NP-hard tree with just 45 vertices by Hell, Ne\v{s}et\v{r}il, and Zhu~\cite{HellNZ96}. The tree $\T$ constructed there is even a \emph{triad}, i.e., a tree with exactly one vertex of degree three and all other vertices of degree one or two.
An NP-hard triad with 39 vertices was found by Barto, Kozik, Mar\'oti, and Niven~\cite{SpecialTriads,SpecialTriadsErratum} using an in-depth analysis of the polymorphisms of triads; they conjectured that their triad is the smallest NP-hard tree (assuming $\operatorname{P}\neq\operatorname{NP}$).
This approach lead to a study of certain classes of trees~\cite{BartoB13,Bulin18}.
Fischer~\cite{Jana} used a computer search and found an NP-hard tree with just 30 vertices (refuting the conjecture of Barto et al.\ mentioned above). Later, independently, Tatarko constructed a 26-vertex NP-hard triad, by manual analysis of polymorphisms. See Table~\ref{tab:history}. 

\renewcommand{\arraystretch}{1.1}
\begin{table}
\begin{center}
\begin{tabular}{llll}
author & year & size & comment \\
\hline
Gutjahr, Welzl, and Woeginger~\cite{GutjahrWW92} & 1992 & 287 & First published\\
Gutjahr~\cite{Gutjahr} & 1991 & 81 & PhD thesis \\
Hell, Ne\v{s}et\v{r}il, and Zhu~\cite{HellNZ96} & 1996 & 45 & Triad \\
Barto, Kozik, Mar\'oti, and Niven~\cite{SpecialTriads} & 2009 & 39 & Triad \\
Fischer~\cite{Jana} & 2015 & 30 & Master thesis \\
Tatarko~\cite{Tatarko} & 2019 & 26 & Triad, Bachelor thesis \\
Present article & 2022 & 22 & Smallest triad \\
Present article & 2022 & 20 & Smallest tree
\end{tabular}
\end{center}
\caption{A time-line of the history of the smallest known NP-hard orientation of a tree.} 
\label{tab:history}
\end{table}

\subsection{Descriptive Complexity} 
Besides the computational complexity of CSPs, the \emph{descriptive complexity} of CSPs has been studied intensively, and leads to a fruitful interplay of finite model theory, graph theory, and universal algebra. Since the results obtained in this context are highly relevant for the open conjectures mentioned above, 
we provide a brief introduction to the  most prominent concepts. A digraph $\HH$ has \emph{tree duality} if for all finite digraphs $\G$, if whenever all trees that map homomorphically to $\G$ also map to $\HH$, then $\G$ maps to $\HH$. It is well-known that a finite digraph $\HH$ has tree duality if and only if the so-called \emph{arc-consistency procedure} solves $\Csp(\HH)$~\cite{FederVardi}. This procedure is of central importance to our work, for many independent reasons that we mention later, and will be introduced in detail in Section~\ref{sect:AC}.

For every finite digraph $\HH$, the arc-consistency procedure for $\Csp(\HH)$ can be formulated as a \emph{Datalog program}~\cite{FederVardi}; Datalog is the fragment of Prolog where function symbols are forbidden. Every Datalog program can be evaluated in polynomial time. Feder and Vardi proved that $\Csp(\HH)$ can be solved by Datalog if and only if $\HH$ has so-called \emph{bounded treewidth duality}; the definition of this concept is similar to the concept of tree duality but we omit it since it is not needed in this article. Bounded treewidth duality can be strengthened to \emph{bounded pathwidth duality}, which corresponds precisely to solvability by a natural fragment of Datalog, namely \emph{linear Datalog}~\cite{Dalmau}. Linear Datalog programs can be evaluated in NL. An even more restricted fragment of Datalog is \emph{linear symmetric Datalog};
such programs can be evaluated in L~\cite{EgriLaroseTessonLogspace}.

A structure $\HH$ has the \emph{ability to count} \cite{AbilityToCount} if $\Csp(\HH)$ can encode, in some natural sense (namely \emph{pp-constructibility} \cite{wonderland,BartoKozikWillardPolymorphisms}), solving systems of linear equations over $\mathbb Z_p$ (for some prime $p$); thus making $\Csp(\HH)$ Mod$_p$L-hard. The ability to count adds substantial complexity to the CSP. Structures that \emph{cannot count} are all tractable and even in Datalog~\cite{Bulatov-BoundedWidth,BartoKozikFOCS09}
and this result, known as the \emph{bounded width theorem}, was an important intermediate step towards the resolution of the Feder-Vardi CSP dichotomy conjecture. Based on this theorem, the lack of the ability to count has a number of equivalent characterizations: \emph{bounded width}, bounded treewidth duality,  definability in Datalog, solvability by \emph{Singleton Arc Consistency}~\cite{Kozik-SLAC}. 

Several important classes of structures exhibit a dichotomy between NP-hardness and the lack of the ability to count (assuming $\operatorname{P}\neq\operatorname{NP}$), which we will refer to as ``easy structures cannot count''. 
Examples include undirected graphs \cite{HellNesetril}, \emph{smooth digraphs} (digraphs without sources and sinks)~\cite{BartoKozikNiven}, \emph{conservative digraphs} (digraphs expanded with all subsets of vertices as unary relations) \cite{HellRafiey-list-homomorphism-digraphs}, \emph{binary conservative structures} (even 3-conservative) \cite{Kazda-binary-conservative}.
We note that this phenomenon also occurs for many large classes of infinite structures $\HH$: for example for all first-order expansions of the basic relations of RCC5~\cite{BodirskyBodorUIP};
see~\cite{Qualitative-Survey} for a survey on the question of which infinite-domain CSPs can be solved in Datalog. In this paper, however, we only consider finite structures. Additionally, for classes of finite structures, 
if the easy structures in the class cannot count, then  the algebraic tractability condition 
for that class can be tested in polynomial time~\cite{MetaChenLarose}.

In \cite{Bulin18} Bul\'in conjectured that easy trees cannot count, establishing this fact for a large yet structurally limited subclass of trees. 
\begin{conjecture}\label{conj:bulin}
Let $\T$ be a tree. If $\T$ has the ability to count, then $\T$ is NP-hard. 
\end{conjecture}

This conjecture (which is rephrasing Conjecture 2 in~\cite{Bulin18}) would answer an open question posed by Hell, Ne\v{s}et\v{r}il, and Zhu~\cite{HNZ} (Open Problem 1 at the end of the article): they asked whether there exists a tractable tree which does not have bounded treewidth duality.

\subsection{Contributions}
\label{sect:results}

In this article, we obtain the following results.

\begin{enumerate}
    \item We find 36 NP-hard trees with 20 vertices; moreover, we prove that all smaller trees and all other trees with 20 vertices are tractable. 
        \item We find four NP-hard triads with 22 vertices, and prove that all smaller triads and all other triads with 22 vertices are tractable.
\item We show that all the trees with at most 20 vertices that are not NP-hard can be solved by Datalog, confirming Conjecture~\ref{conj:bulin} for trees with at most 20 vertices. 
    \item We find a tree with 19 vertices that can not be solved by arc consistency, and prove that all smaller trees and all other trees with 19 vertices can be solved by arc consistency.
    \item We find 8 NL-hard trees with 12 vertices; moreover, we prove that all smaller trees and all other trees with 12 vertices are in L. 
\end{enumerate}

Even though we draw from the results of the universal-algebraic approach to the CSP which led to the theorems of Bulatov and of Zhuk, and use state-of-the-art computers for our computations, these tasks remain challenging due to the huge number of trees: for example, even considered up to isomorphism, there are
139,354,922,608 trees with 20 vertices (see Table~\ref{table:otrees}), which is prohibitive even if we could 
test the algebraic tractability condition
within milliseconds. 
Several further contributions of this article are related to the way we managed to overcome these difficulties.

A well-known key simplification is to only consider trees that are \emph{cores}; a digraph $\HH$ is a core if every homomorphism from $\HH$ to $\HH$ is injective. Two graphs $\HH_1$ and $\HH_2$ are called \emph{homomorphically equivalent} if there is a homomorphism from $\HH_1$ to $\HH_2$ and vice versa. 
Clearly, in that case $\Csp(\HH_1)=\Csp(\HH_2)$ and so $\HH_1,\HH_2$ are either both  tractable or both NP-hard. 
It is easy to see that every finite digraph $\HH$ is homomorphically equivalent to a core digraph $\HH'$, which is unique up to isomorphism. Moreover, if $\HH$ is a tree, then $\HH'$ is a tree as well (and its size is smaller or equal to the size of $\HH$). Hence, it suffices to work with trees that are cores.
However, Hell and Ne\v{s}et\v{r}il proved that deciding whether a given digraph is a core is coNP-complete~\cite{Cores}.
Remarkably, the following fact seems to be unnoticed in the literature.

\begin{enumerate}[resume]
    \item There is a polynomial-time algorithm for deciding whether a tree is a core. 
\end{enumerate}

There are far too many trees with at most 20 vertices to run the core test on each of them. Our next contribution is a method to generate the core trees more directly, rather than generating all trees and then discarding the non-cores (the details can be found in Section~\ref{sect:gen}). 
Applying our method we were able to construct all trees that are cores up to size 20, and all triads that are cores up to size 22, which was essential to achieve the results 1.-5.\ above.

\begin{enumerate}[resume]
    \item 
    We computed the number of trees that are cores for sizes up to 20 (see Table~\ref{table:otrees}). In particular, there are 779268 core trees of size 20.
\end{enumerate}
These are still too many to be tested for the algebraic tractability condition 
if this is implemented naively. We therefore use results from the universal algebraic approach to first run more efficient tests for certain sufficient conditions,
such as the existence of a binary symmetric polymorphism, and only run the full test for Kearnes-Markovi\'{c}-McKenzie polymorphisms if the simpler conditions all fail; this will be explained in Section~\ref{sect:specific-polymorphisms}.

Finally, we identify trees that are important `test-cases' for the open problems that have been mentioned earlier. 
\begin{enumerate}[resume]
    \item We computed the two smallest trees that are not known to be in NL; they have 16 vertices, and they are the smallest trees that do not have a majority polymorphism. 
    \item We computed 28 smallest trees that are candidates for failing the condition in Conjecture~\ref{conj:NL} and hence 
    might be P-hard (and that are thus candidates for not being in NL, unless $\operatorname{NL}=\operatorname{P}$); they have 18 vertices. 
\end{enumerate}

\subsection{Outline of the Article}
Basic notation and terminology about directed and undirected graphs and homomorphisms is introduced in Section~\ref{sect:prelims}. This section also presents a brief description of the arc-consistency procedure which plays an important role in several of our results. 
In Section~\ref{sect:cores} we explain how to use the arc-consistency procedure to efficiently test whether a given tree is a core. 
In Section~\ref{sect:gen} we present our method to generate all trees that are cores (directly, without having to discard too many non-cores in the process).  
We then use these trees to make extensive experiments about their computational and descriptive complexity.
For this, we need to introduce important polymorphism conditions and related facts from universal algebra (Sections~\ref{sect:ua} and \ref{sect:specific-polymorphisms}). Finally, the results of our experiments as announced in Section~\ref{sect:results} can be found in Section~\ref{sect:exp}. 

\section{Graphs,  Digraphs, Homomorphisms}
\label{sect:prelims}
For the definition of \emph{relational structure} we refer to any text-book in mathematical logic; note that we allow the signature of structures to be infinite (but the constraint satisfaction problem is only defined for relational structures with finite relational signature). Since we work most of the time with digraphs, we present the basic definitions only for digraphs; most of them generalize to relational structures in a straightforward way. We use standard terminology for graphs and undirected graphs as introduced e.g.\ in~\cite{Diestel}. All graphs we consider are finite.
A digraph is a pair $\HH = (H;E)$ where $H$ is a nonempty set and $E=\edgesH \subseteq H^2$ is a set of (directed) edges.  A~\emph{(simple, undirected) graph} is a pair $\HH = (H;E)$ where $H$ is a nonempty set and $E=\edgesH \subseteq {H \choose 2}$ is a set of two-element subsets of $H$.
An \emph{orientation of} a graph $\G$ is a digraph $\digraph O$ such that $O=G$, $(x,y) \in \edges{\digraph O}$ implies
$\{x,y\} \in \edges{\G}$, 
and for every $\{x,y\} \in \edges{\G}$ either $(x,y) \in \edges{\digraph O}$ or $(y,x) \in \edges{\digraph O}$, but not both. If $\HH$ is a digraph, then the \emph{reverse} of $\HH$ is the digraph $\HH^R = (H;E^R)$ where $E^R = \{(y,x) \mid (x,y) \in E\}$. 
The operation that obtains $\HH^R$ from $\HH$ is called \emph{edge reversal}.

If $\G$ and $\HH$ are digraphs, then a \emph{homomorphism} from $\G$ to $\HH$ is a map $h \colon G \to H$ such that for all $(x,y) \in \edges{\G}$ we have $(h(x),h(y)) \in \edges{\HH}$. We write $\Csp(\HH)$ (for \emph{constraint satisfaction problem})  
for the class of all finite digraphs $\G$ which admit a homomorphism to $\HH$. A homomorphism from $\HH$ to $\HH$ is called an \emph{endomorphism of $\HH$}. 
A finite digraph $\HH$ is called a \emph{core} if all endomorphisms of $\HH$ are injective. It is easy to see that an injective endomorphism of $\HH$ must 
in fact be bijective and an \emph{automorphism}, i.e., an isomorphism between $\HH$ and $\HH$. It is easy to see that every finite digraph $\G$ is homomorphically equivalent to a finite core digraph $\HH$, and that this core digraph is unique up to isomorphism~\cite{HNBook}, hence $\HH$ will be called \emph{the} core of $\G$.

An \emph{undirected tree} is a connected  undirected graph without cycles. If $u, v \in T$ and $\T$ is an undirected tree, then  there exists a unique path $\digraph P$ from $u$ to $v$ in $\T$;
the number of edges of $\digraph P$ is denoted by $\dist(u,v)$. 
A vertex $v\in T$ is called a \emph{center} of $\T$ if $v$ lies in the middle of a longest path in $\T$.
An edge $e\in \edges{\T}$ is called a \emph{bicenter} of $\T$ if $e$ is the middle edge of a longest path in $\T$. 
We will use the following classical result.

\begin{theorem}[Jordan (1869)]\label{thm:jordan} 
An undirected tree $\T$ has exactly one center or one bicenter.
\end{theorem}

If $\digraph O$ is an orientation of a tree and $u,v \in O$,  then $\dist(u,v)$ (center of $\digraph O$,  bicenter $\digraph O$)
is meant with respect to the underlying undirected tree. 
As mentioned in the introduction in this article an orientation of a finite tree will simply be called a \emph{tree}. 
A digraph $\HH$ is \emph{balanced} if its vertices can be organized into \emph{levels}, that is, there exists a function $\levelfunction\colon H\to\mathbb N$ such that $\level v=\level u + 1$ for all $(u,v)\in\edgesH$ and 
the smallest level is 0. The \emph{height} of $\HH$ is the maximum level. Note that trees are balanced and observe that if $\G$ and $\HH$ are balanced of the same height, and $\G$ is connected, then any homomorphism from $\G$ to $\HH$ must preserve levels, that is, $\level{v}=\level{h(v)}$ for all $v\in G$.
A \emph{rooted tree} is a tuple $(\T,r)$, where  $\T$ is a tree and $r\in T$; $r$ is then called the \emph{root} of $\T$. 
A rooted tree $(\T,r)$ is called a \emph{rooted core} if every endomorphism of $\T$ that fixes $r$ is injective. 
The \emph{depth} of a rooted tree $(\T,r)$ is $\max\{\operatorname{dist}(r,v)\mid v \in T\}$.

\subsection{The Arc-consistency Procedure}
\label{sect:AC}
One of the most efficient algorithms employed by constraint solvers to reduce the search space 
is the \emph{arc-consistency} procedure.
In the graph homomorphism literature, the algorithm is sometimes
called the \emph{consistency check algorithm}.
The arc-consistency procedure is important for us for several reasons:
\begin{itemize}
    \item It plays a crucial role for efficiently deciding whether a given tree is a core (Section~\ref{sect:cores}). 
    \item It is well suited for combination with exhaustive search to prune the search space, and this will be relevant in Section~\ref{sect:ua}.
    \item It is an important fragment of Datalog of independent interest from the point of view of the CSP theory (see Section~\ref{sect:AC-Pol}), and we will later perform experiments to compute the smallest tree that cannot be solved by arc consistency (Section~\ref{sect:NAC}). 
\end{itemize}
We need to give a short description of the procedure. 

Let $\G$ and $\HH$ be finite digraphs. We would like to determine whether there exists a homomorphism from $\G$ to $\HH$. 
The idea of the arc-consistency procedure is to maintain for each $x \in G$ a set $L(x) \subseteq H$.
Informally, each element of $L(x)$ represents a candidate for an image of $x$ under a homomorphism from $\G$ to $\HH$.
The algorithm initializes each list $L(x)$ with $H$ and successively removes vertices from these lists; it only removes 
a vertex $u \in H$ from $L(x)$ if there is no homomorphism from $\G$ to $\HH$ that maps $x$ to $u$. 
To detect vertices $x,u$ such that $u$ can be removed from $L(x)$, the algorithm uses two rules (in fact, one rule and
a symmetric version of the same rule): if $(x,y) \in \edges{\G}$, then
\begin{align*}
\text{remove $u$ from $L(x)$ if there is no $v \in L(y)$
with $(u,v) \in \edges{\HH}$;}\\
\text{remove $v$ from $L(y)$ if there is no $u \in L(x)$
with $(u,v) \in \edges{\HH}$.}
\end{align*}
If eventually we cannot remove any vertex from any list with these rules any more, the digraph $\G$ together with 
the lists for each vertex is called \emph{arc-consistent}.
Note that formally we may view $L$ as a function $L \colon G \to 2^H$  from the vertices of $\G$ to sets of vertices of $\HH$.

Note that we may run the algorithm also on digraphs $\G$ where for some $x \in G$ the list $L(x)$ is already set to some subset of $H$. 
In this setting, the input consists of $\G$ and the given lists, and we are looking for a homomorphism $h$ from $\G$ to $\HH$ such that 
for every $x \in G$ we have $h(x) \in L(x)$. 
The pseudocode of the entire arc-consistency procedure is displayed in Algorithm~\ref{alg:ac}. The standard 
arc-consistency procedure $\AC_\HH(\G)$ is then obtained by calling $\AC_\HH(\G,L)$ with $L(x) := H$ for all $x \in G$.

\RestyleAlgo{ruled}
\SetAlgoVlined{}
\begin{algorithm}
\DontPrintSemicolon{}
\SetKwInOut{Input}{input}\SetKwInOut{Data}{data}
\SetKw{Remove}{remove}
\Input{a finite digraph $\G$}
\Data{a list $L(x) \subseteq H$ for each vertex $x \in G$}
\caption{$\AC_\HH(\G,L)$ }\label{alg:ac}
\Repeat{no list changes}{
  \ForEach{$(x,y) \in \edges{\G}$}{
    \If{there is no $v \in L(y)$ with $(u,v) \in \edges{\HH}$}{\Remove{$u$ from $L(x)$}}
    \If{there is no $u \in L(x)$ with $(u,v) \in \edges{\HH}$}{\Remove{$v$ from $L(y)$}}
  }
  \If{$L(x)$ is empty for some vertex $x \in G$}{\bf reject}
}
\end{algorithm}

Clearly, if the algorithm removes all vertices from one of the lists, then there is no homomorphism from $\G$ to $\HH$. 
It follows that if $\AC_\HH$ rejects $\G$, then there is no homomorphism from $\G$ to $\HH$.  
The converse implication does not hold in general.
For instance, let $\HH$ be the loopless digraph with two vertices and two edges, denoted $\digraph K_2$, and let $\G$ 
be $\digraph K_3 = (\{0,1,2\};\neq)$. In this case, $\AC_\HH$ does not remove any vertex from any list, but obviously 
there is no homomorphism from $\digraph K_3$ to $\digraph K_2$.

The arc-consistency procedure can be implemented so that it runs in $O(|\edges{\G}| \cdot |H|^3)$, e.g.~by Mackworth's AC-3 algorithm~\cite{Mackworth}.

\section{Cores of Trees}
\label{sect:cores}
Recall that the problem of deciding whether a given digraph is a core is coNP-complete~\cite{Cores}. 
The following theorem implies that 
whether a given  finite tree is a core can be tested in polynomial time.

\begin{theorem}\label{thm:ac} 
Let $\T$ be a finite tree. Then the following are equivalent.
\begin{enumerate}
    \item $\T$ is a core;
\item $\End(\T) = \{\id_{T}\}$;
\item $\AC_\T(\T)$
terminates such that the list for each vertex contains a single element.
\end{enumerate}
\end{theorem}

\begin{corollary}\label{cor:treecoreP}
There is a polynomial-time algorithm to decide whether a given finite tree is a core. 
\end{corollary}

We first prove the following two useful lemmata.

\begin{lemma}\label{lem:ac-new}
Let $\T$ be a finite tree and let $\HH$ be a finite digraph such that $\AC_\HH(\T)$ does not reject. Let $t\in T$, and 
let $a\in H$ be such that $a\in L(t)$ after running $\AC_\HH(\T)$. Then there is a homomorphism $h\colon \T\to \HH$ such that $h(t)=a$. 
\end{lemma}

\begin{proof}
Let $S$ be a maximal subtree of $\T$ such that $t\in S$ and there exists a partial homomorphism $h\colon S\to\HH$ with $h(t)=a$. 
If $S\neq T$, then there exists $x\in S$ and $y\in T\setminus S$ such that either $(x,y)$ or $(y,x)$ is an edge in $\T$; without 
loss of generality assume that $(x,y)\in\edges{\T}$. Because the value $u=h(x)$ was not removed from $L(x)$ when running $\AC_\HH(\T)$, 
it follows that there exists $v\in L(y)$ such that $(u,v)\in\edges{\digraph H}$. But then setting $h(y)=v$ extends $h$ to a 
partial homomorphism from $S\cup\{y\}$ to $\HH$ contradicting maximality of the subtree $S$.
\end{proof}

\begin{lemma}\label{lem:rootedTreeAutoImpliesEndo}
Let $(\digraph T,r)$ be a rooted tree  with an automorphism that is not the identity. Then $(\T,r)$ has a non-injective endomorphism.
\end{lemma}
\begin{proof}
We prove the statement by induction on the number of vertices of $\digraph T$.
Consider the components of the graph obtained from $\T$ by deleting $r$. If there is a component $C$ such that $h$ does not map $C$ into itself, then the mapping which agrees with $h$ on $C$ and which fixes all other vertices of $\T$ is a non-injective endomorphism of $(\T,r)$.

If each component $C$ is mapped by $h$ into itself, then each $h|_C$ is an automorphism of $(C,r_C)$, where $r_C$ is the unique neighbor of $r$ that lies in $C$. Since $h$ is not $\operatorname{id}_T$ there must be some $C$ such that $h|_C$ is not $\operatorname{id}_C$ and by induction hypothesis there exists a non-injective endomorphism $h'$ of $(C,r_C)$. Since $h'(r_C)=r_C$ the mapping which extends $h'$ to $T$ by fixing all other vertices of $\T$ is a non-injective endomorphism of $(\T,r)$.
\end{proof}

\begin{proof}[Proof of Theorem~\ref{thm:ac}]
We prove the equivalence of 1.\ and 2., and then the equivalence of 2.\ and 3.

Clearly, 2.\ implies 1. Conversely, suppose that $\T$ has an endomorphism $h$ which is not the identity map. If $h$ is not injective, then $\digraph T$ is not a core and we are done. 
Hence, suppose that $h$ is an automorphism. Note that by Theorem~\ref{thm:jordan}, if $\T$ has a center, then $h(c)=c$ and if $\T$ has a bicenter $(x,y) \in E(\T)$, then $h(\{x,y\})=\{x,y\}$. In the latter case, since $(y,x)\notin E(\T)$, we must have $h(x)=x$ and $h(y)=y$.
In both cases $h$ has a fixed point $r$ and $h$ is an automorphism of $(\T,r)$. By Lemma~\ref{lem:rootedTreeAutoImpliesEndo}, $\T$ has a non-injective endomorphism and is therefore not a core.

To prove that 2.\ implies 3., we prove the contrapositive. 
Suppose that $\AC_\T(\T)$ terminates with $|L(x)|>1$ for some $x \in T$. Then there exists $y\in L(x),y\neq x$. Lemma~\ref{lem:ac-new} implies that there is an endomorphism $h$ of $\T$ such that $h(x) = y$ and thus $h\neq\id_{T}$.

To see that 3.\ implies 2., note that $L(x)=\{x\}$ since because of the identity endomorphism $\id_{T}$, $x$ cannot be removed from $L(x)$. Therefore, for any endomorphism $h\colon\T\to\T$ and any $x\in T$ we must have $h(x)\in L(x)$, and so $h(x)=x$ and $h=\id_{T}$.
\end{proof}

\section{Generating all Core Trees}
\label{sect:gen} 
In this section we present an algorithm to generate all core trees with $n$ vertices up to isomorphism. 
To this end, we first present a known algorithm that generates all trees with $n$ vertices up to isomorphism \cite{103287}.
Later we explain how to modify this algorithm to directly generate core trees.

We also refer to the isomorphism classes of trees as \emph{unlabeled trees}, as opposed to \emph{labeled trees}, which are trees with vertex set $\{1,\dots,n\}$ for some $n \in {\mathbb N}$. 
The difference between the enumeration of labelled and  unlabeled trees is significant:
while the number of labelled trees is Sloane's integer sequence A097629, given by $2 (2n)^{n-2}$, 
the number of unlabeled trees  is Sloane's integer sequence A000238, which grows
asymptotically as $c d^n / n^{5/2}$
where $d \approx 5.6465$ and $c \approx 0.2257$ are constants;
the initial terms are shown in Table~\ref{table:otrees}. 
However, these numbers are still too large to apply the core test to all the unlabeled trees separately. 
The number of unlabeled trees that are cores is again much smaller. 
We therefore present a modification of the generation algorithms that allows us to generate  unlabeled trees that are cores directly without enumerating all unlabeled trees.

Let $\geq$ be some total order on all rooted trees that linearly extends the order by depth. 
The idea of the algorithm is to generate all unlabeled rooted  trees with at most $n-1$ vertices and then use Theorem~\ref{thm:jordan}.

\begin{algorithm}[ht]
\SetKwInOut{Input}{input}\SetKwInOut{Output}{output}
\caption{$\operatorname{GenerateTrees}$}\label{alg:genTrees}
\DontPrintSemicolon{}
\Input{a positive integer $n$}
\Output{a list of trees with $n$ vertices}
\Begin{
$\mathrm{Trees} \leftarrow \emptyset$\;
\tcp{bicenter}
\ForEach{$(\T_1,r_1),(\T_2,r_2)$ rooted trees where
  $|T_1| + |T_2| = n$ and $\operatorname{depth}(\T_1,r_1)=\operatorname{depth}(\T_2,r_2)$
}{
  $T := T_1\uplus T_2$\\
  $E := \{(r_1,r_2)\}\uplus \edges{\T_1}\uplus \edges{\T_2}$\\
  $\T := (T;E)$\\ 
  $\mathrm{Trees} \leftarrow \mathrm{Trees} \cup \{\T\}$
}

\hspace{0.5cm}

\tcp{center}
\ForEach{ $(\T_1,r_1)\geq(\T_2,r_2)\geq\dots\geq(\T_m,r_m)$ rooted trees where
$|T_1|+\dots+|T_m|=n-1$ and $\operatorname{depth}(\T_1,r_1)=\operatorname{depth}(\T_2,r_2)$
}{
  \ForEach{ $s\in\{0,1\}^m$ where $i<j$ and $(\T_i,r_i)=(\T_j,r_j)$ imply $s_i\leq s_j$}{
  $T := \{r\}\uplus T_1\uplus\dots\uplus T_m$\\
  $E := \{(r_i,r)\mid r_i=1\}\uplus\{(r,r_i)\mid r_i=0\}\uplus E(\T_1)\uplus\dots\uplus E(\T_m)$\\
  $\T := (T;E)$\\
  $\mathrm{Trees} \leftarrow \mathrm{Trees} \cup \{\T\}$
 }
}
\Return{$\mathrm{Trees}$}
}
\end{algorithm}

\begin{algorithm}
\caption{$\operatorname{GenerateRootedTrees}$}\label{alg:genRootedTrees}
\SetKwInOut{Input}{input}\SetKwInOut{Output}{output}
\Input{two positive integers $n$, $d$}
\Output{a list of rooted trees with $n$ vertices and depth
$d$}
\Begin{
\If{$n=0$}{
\Return{$\emptyset$}}
\If{$n=1$}{
\Return{$\{((\{r\};\emptyset),r)\}$}}
$\mathrm{RootedTrees} \leftarrow \emptyset$\\
\ForEach{$(\T_1,r_1)\geq\dots\geq(\T_m,r_m)$ rooted trees where
 $|T_1|+\dots+|T_m|=n-1$  and $\operatorname{depth}(\T_1,r_1)=d-1$}{
    \ForEach{$s\in\{0,1\}^m$ where $i<j$ and $(\T_i,r_i)=(\T_j,r_j)$ imply $s_i\leq s_j$}{
  $T := \{r\}\uplus T_1\uplus\dots\uplus T_m$\\
  $E := \{(r_i,r)\mid r_i=1\}\uplus\{(r,r_i)\mid r_i=0\}\uplus \edges{\T_1}\uplus\dots\uplus \edges{\T_m}$\\
  $\T := (T;E)$\\
  $\mathrm{RootedTrees} \leftarrow \mathrm{RootedTrees} \cup \{(\T,r)\}$
  }
    }
\Return{$\mathrm{RootedTrees}$}
}
\end{algorithm}

It is easy to verify that Algorithm~\ref{alg:genRootedTrees} produces all unlabeled rooted trees with $n$ vertices and depth~$d$.
Analogously,  Algorithm~\ref{alg:genTrees} generates all unlabeled trees with $n$ vertices.
Remarkably, there are no isomorphism checks necessary and Algorithm~\ref{alg:genTrees} runs in linear time in the number of unlabeled trees with $n$ vertices plus the number of unlabeled rooted trees with at most $n-1$ vertices. 

Let us make some observations.
Let $(\T_1,r_1)\geq\dots\geq(\T_m,r_m)$ be rooted trees, $s\in\{0,1\}^m$, $T := \{r\}\uplus T_1\uplus\dots\uplus T_m$, and $E := \{(r_i,r)\mid s_i=1\}\uplus\{(r,r_i)\mid s_i=0\}\uplus E(\T_1)\uplus\dots\uplus E(\T_m)$.
\begin{itemize}
    \item A rooted tree $(\T,r)$ is a rooted core if and only if $\AC_\T(\T,L)$, where $L(r)=\{r\}$ and $L(x)=T$ for $x\in T\setminus\{r\}$, terminates such that the list for each vertex contains a single element.
    \item By Corollary \ref{cor:treecoreP}, testing whether a (rooted) tree is a (rooted) core can be checked in polynomial time using the arc-consistency procedure. 
    \item If $(\T,r)$ is a rooted core, then  $(\T_i,r_i)$ is a rooted core for every $i$.
    \item If $\T$ is a core and $r$ is its center, then  $(\T_i,r_i)$ is a rooted core for every $i$.
    \item If $(T_1\uplus T_2;\{(r_1,r_2)\}\uplus \edges{\T_1}\uplus \edges{\T_2})$ is a core and $(r_1,r_2)$ is its bicenter, then $(\T_1,r_1)$ and $(\T_2,r_2)$ are rooted cores.
    \item If two trees that are cores are homomorphically equivalent, then they are isomorphic. 
\end{itemize}

To generate all oriented trees that are cores we slightly modify both algorithms. In both functions we only add trees to the output if they are cores or rooted cores, respectively.
By the above observations, these modified algorithms generate each tree with $n$ vertices that is a core exactly once. We do now know whether our algorithm is a polynomial-delay generation procedure for unlabeled core trees. In practice, it is fast enough to generate all core trees with at most 20 vertices within  reasonable time (see Section~\ref{sect:exp}).

\section{Polymorphism Conditions}
\label{sect:ua}
In this section we introduce basic facts from the universal-algebraic approach that are essential for obtaining our results. 
If $\HH$ is a digraph and $k \geq 1$ is an integer, then $\HH^k$ denotes the \emph{$k$-th categorical power of $\HH$}, i.e., the digraph with vertex set $H^k$ and edge set 
$$\{((u_1,\dots,u_k),(v_1,\dots,v_k)) \mid (u_1,v_1),\dots,(u_k,v_k) \in \edges{\HH} \}.$$ 
A \emph{polymorphism of $\HH$} is a homomorphism from $\HH^k$
to $\HH$, for some $k \geq 1$. Clearly, every projection, i.e., every operation of the form $(x_1,\dots,x_k) \mapsto x_i$, for some fixed $i \leq k$, is a polymorphism for every digraph $\HH$.
Of particular interest to us will be polymorphisms that satisfy certain sets of identities, introduced in Section~\ref{sect:height-one}. 
The connection between such identities and computational complexity is described in Section~\ref{sect:pp}.

\subsection{Primitive Positive Constructions} 
\label{sect:pp} 
Primitive positive definitions are a natural type of gadget construction which can be used to 
obtain logspace reductions between CSPs. 
Central to the algebraic theory of the CSP is the fact that the possibility of such an encoding can be determined from the polymorphisms of the respective templates. 
Let $\structure A = (A;R_1,\dots,R_n)$ be a relational structure. A relation $S\subseteq A^n$ is \emph{primitive positive (pp-) definable} from $\structure A$ if it can be defined (without parameters) by a first order formula which only uses the predicate symbols $R_1,\dots,R_n$, the equality predicate, conjunction, and existential quantification. 
 
 \emph{Primitive positive constructions} are a more powerful generalization of primitive positive definitions. 
A relational structure $\structure B=(B;S_1,\dots,S_m)$ is \emph{pp-constructible} from $\structure A$  if there exists $k>0$ and a structure $\structure B'=(B';S'_1,\dots,S'_m)$ which is  homomorphically equivalent to $\structure B$ such that $B'=A^k$ and for every $i\in\{1,\dots,m\}$ the relation $S_i\subseteq {B'}^{r_i}$, when considered as a $kr_i$-ary relation on $A$, is pp-definable from $\structure A$. 
In this case there is a logspace reduction from $\Csp(\structure B)$ to $\Csp(\structure A)$.

If $\structure A$ and $\structure B$ are finite structures, then $\structure B$ is pp-constructible from $\structure A$ if and only if $\structure B$ satisfies every \emph{height-one condition satisfied by $\structure A$}; these concepts will be introduced in the next section.  
This algebraic characterization of pp-constructibility was shown in~\cite[Theorem 1.3, Corollary 4.7]{wonderland} (see also \cite[Theorem 38, Corollary 20]{BartoKozikWillardPolymorphisms}).

\subsection{Linear Conditions and Height-one Conditions}
\label{sect:height-one}
Height-one conditions and linear conditions are particular types of \emph{strong Maltsev condition}~\cite{KearnesKiss} that are essential for the algebraic approach to CSPs. 
If $f$ is a function symbol of arity $k$ and $h$ is a function symbol of arity $\ell$, and $\sigma \colon \{1,\dots,k\} \to \{1,\dots,n\}$ and 
$\rho \colon \{1,\dots,\ell\} \to \{1,\dots,n\}$
are functions, 
then an expression of the form
$$f(x_{\sigma(1)},\dots,x_{\sigma(k)}) \approx g(x_{\rho(1)},\dots,x_{\rho(\ell)})$$
is called a \emph{height-one identity}. For the purposes of this paper, a finite set of height-one identities will be called a \emph{height-one condition}. Height-one conditions are important because of the mentioned tight link with pp-constructibility (Section~\ref{sect:pp}). More generally, an identity is \emph{linear} if each side has one or zero occurrences of function symbols, i.e., it is either height-one, or of the form $f(x_{\sigma(1)},\dots,x_{\sigma(k)})\approx x_j$, $x_i\approx g(x_{\rho(1)},\dots,x_{\rho(\ell)})$, or $x_i\approx x_j$.
Then a \emph{linear condition} is a finite set of linear identities.

A set of operations $\mathcal F$ on some domain $D$ \emph{satisfies} a linear condition $\Sigma$ if we can interpret every function symbol appearing in $\Sigma$ as an operation from $\mathcal F$ so that for every identity in $\Sigma$, the left-hand side of the identity and the right-hand side of the identity evaluate to the same element under all possible substitutions of variables by elements of $D$.
We say that a structure $\structure A$ \emph{satisfies} a linear condition if the set of all polymorphisms of $\structure A$ satisfies it.

An example of a height-one condition consisting of a single height-one identity involving a single binary function symbol $f$ is 
$$ \{f(x_1,x_2) \approx f(x_2,x_1) \}$$
which is satisfied by $\structure A$ if and only if $\structure A$ has a binary symmetric polymorphism. Since $x_1,x_2,x_3,\dots$ are just variable names we sometimes use $x,y,z$, etc.\ instead.

An operation $f \colon D^k \to D$ is called \emph{idempotent} if it satisfies the identity  $f(x,\dots,x) \approx x$. Note that this identity is linear but not height-one.

\begin{remark}
\label{rem:idempotence} 
It is well-known and easy to see that the polymorphisms of a finite core digraph satisfy a height-one condition 
if and only if its idempotent polymorphisms satisfy the condition. 
\end{remark}

Linear conditions have been introduced and studied first, in particular motivated by the fact that idempotence of operations can only be expressed in linear conditions, but not in height one conditions, and that idempotence plays a central role in many classical areas of universal algebra. 
Also note that in the setting of finite-domain CSPs we may use the idempotence assumption for free because of Remark~\ref{rem:idempotence} and the fact that the core of a structure has the same CSP and satisfies the same height-one identities (which is not true for linear conditions in general).

\subsection{The Indicator Construction}
\label{subsect:testing-linear}

The question whether a given digraph $\HH$ has polymorphisms that satisfy a given height-one or even linear condition can be tested algorithmically as follows. To illustrate the well-known idea, suppose that the given set consists of a single identity, namely $f(x,y) \approx f(y,x)$. 
We then compute  $\HH^2$ and contract every vertex of $\HH^2$ of the form $(x,y)$ with the vertex $(y,x)$. The resulting digraph $\indicator{\HH}$ will be called the \emph{indicator digraph} for the height-one condition. 
We finally search for a homomorphism from $\indicator{\HH}$ to $\HH$. Note that $\indicator{\HH}$ may be viewed as an instance of $\Csp(\HH)$, and that the homomorphisms from $\indicator{\HH}$ to $\HH$ are in 1-1 correspondence with the binary symmetric polymorphisms of $\HH$ (see Definition \ref{def:wnu}).

Analogously we may proceed for any other height-one condition: to compute $\indicator{\HH}$, we construct for each function symbol the categorical power of $\HH$ of the corresponding arity, take their disjoint union, and then identify vertices as dictated by the identities. 
Clearly, the size of the indicator digraph grows exponentially with the arity of the function symbols in the condition and linearly with number of function symbols in the condition so we generally prefer conditions where the function symbols are of low arity even if the number of function symbols is large. 

Linear conditions can be tested in the following way. Note that the left- and right-hand sides of identities can be switched, and that identities of the form $x_i\approx x_j$ are only satisfied in one-element structures. Therefore, we may assume that every identity is either height-one or of the form $f(x_{\sigma(1)},\dots,x_{\sigma(k)})\approx x_j$. First, construct the indicator digraph $\indicator{\HH}$ using only the height-one identities when identifying vertices. Then, for every identity that is not height-one,  find every vertex of $\indicator{\HH}$ that comes from a tuple of vertices of $\HH$ matching the left-hand side and set its value to the vertex of $\HH$ given by the right-hand side. For example, if the identity is $f(x,y,x)\approx x$, we require that $(a,b,a)\in\indicator{H}$ must be mapped to $a$, for every $a,b\in H$. In this way, we obtain an instance of the \emph{$\HH$-precoloring extension problem} (see also the discussion in Section \ref{sect:discussion}). It is well known that for cores, this problem is logspace-equivalent to $\Csp(\HH)$ \cite{JBK,wonderland}. Moreover, it is particularly easy to implement within the arc-consistency procedure, see the next section.
For balanced digraphs, another important improvement---based on the decomposition into levels---that can be applied for many height-one conditions is described in Section~\ref{sect:level-wise} below.

\subsection{Arc Consistency with Exhaustive Search}
Chen and Larose~\cite{MetaChenLarose} give a polynomial-time procedure to test for the existence of polymorphisms satisfying a given linear condition, providing that the linear condition implies the algebraic tractability condition. Their procedure uses the algorithm of Zhuk which is known to be \emph{uniform} in the sense that it runs 
in polynomial time even if we assume that also the template $\HH$ is part of the input to the algorithm for $\Csp(\HH)$~\cite{ZhukFVConjecture}.

However, in our approach we avoid implementing and running the algorithm of Zhuk (which is complex and whose running time is a polynomial of a yet unknown degree). Instead, we run the arc-consistency procedure for $\HH$ on the indicator digraph $\indicator{\HH}$ and then perform an exhaustive search. While this procedure is not (provably) in P, it is very efficient in practice.

We initialize the lists for vertices of $\indicator{\HH}$ with preset values dictated by non-height-one identities, as explained above. Additionally, for every $u \in H$, we initialize the list for every vertex of $\indicator{\HH}$ of the form $(u,\dots,u)$ with $\{u\}$ (since it suffices to look for idempotent polymorphisms as we have explained in Remark~\ref{rem:idempotence} and this reduces the search space). For the remaining vertices of $\indicator{\HH}$, the lists are initialized to $H$.

If $\AC_\HH$ detects an inconsistency, we can be sure that no polymorphisms satisfying the linear condition exist. Otherwise, we select some vertex $x \in \indicator{H}$, and set $L(x)$ to ${u}$ for some $u \in L(x)$. 
Then we proceed recursively with the resulting lists. If $\AC_\HH$ now detects an empty list, we backtrack, but remove $u$ from $L(x)$. Finally, if the algorithm does not detect an empty list at the first level of the recursion, we end up with singleton lists for each vertex $x \in \indicator{H}$, which defines a homomorphism from $\indicator{\HH}$ to $\HH$. This homomorphism can then be interpreted as polymorphisms satisfying the linear condition.

There are numerous methods to speed up this backtracking procedure. One of the best known is called 
Maintaining Arc Consistency (MAC) \cite{Sabin1994ContradictingCW}. This family of algorithms has the arc-consistency procedure
at its core and takes advantage of the incremental design of the backtracking procedure by maintaining data structures which help to 
reduce the number of consistency checks. Another common way to speed up the search procedure is to choose the vertex $x \in \indicator{\HH}$ 
that has a list of smallest size.

\subsection{Level-wise satisfiability}\label{sect:level-wise}

If $\HH$ is a balanced digraph (in particular, a tree), the test from the previous section can sometimes be significantly simplified. We say that a linear condition is \emph{level-wise satisfied} if we can interpret the function symbols as polymorphisms of $\HH$ in such a way that for every level in $\HH$, the identities are satisfied under all evaluations of variables by vertices from that level.

When testing whether a linear condition is level-wise satisfied, we do not need to construct the full indicator digraph. Instead, for every function symbol (say of arity $k$) we construct only the subgraph of $\HH^k$ consisting of $k$-tuples of same-level vertices. Note that this is a union of connected components of $\HH^k$ and that polymorphisms can be defined the first projection on the remaining connected components of $\HH^k$.

While we do not have a general construction, for many linear conditions relevant to the complexity of the CSP we can show that if a linear condition is level-wise satisfied in $\HH$, then it is satisfied in $\HH$. The idea is to start with polymorphisms satisfying the identities level-wise, and then redefine those polymorphisms for tuples of vertices that are not all on the same level, in such a way as to satisfy the identities. We will introduce several such concrete constructions in the next section. Similar constructions have appeared in \cite{SpecialTriads,BartoB13,Bulin18,BulinDelicJacksonNiven}.
This optimization is particularly useful when testing the condition $\TS(n)$ for all $n$; see Section~\ref{sect:NAC}.

\section{Specific Polymorphism Conditions}\label{sect:specific-polymorphisms}
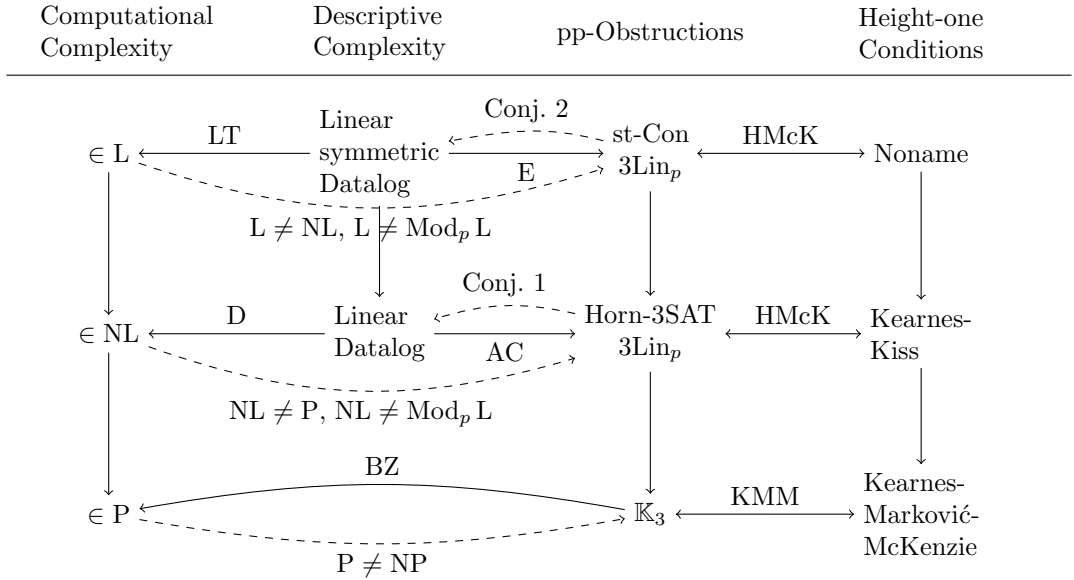
\begin{figure}
    \centering
    \begin{tikzpicture}[scale=0.8]
        \node[align=left] at (0.3,0) {Computational\\ Complexity};
        \node[align=left] at (4.5,0) {Descriptive\\ Complexity};
        \node[align=left] at (9,0) {pp-Obstructions};
        \node[align=left] at (13.5,0) {Height-one\\ Conditions};
        
        \node (L1) at (0,-2) {$\in \operatorname{L}$};
        \node[align=left] (LSD1) at (4.5,-2) {Linear\\symmetric\\Datalog};
        \node (Leq1) at (9,-2) {\parbox{1cm}{\centering $\stcon$ $\TLinP$}};
        \node[align=left] (HM1) at (13.5,-2) {Noname};
        
        \node (L2) at (0,-5) {$\in \operatorname{NL}$};
        \node[align=left] (LSD2) at (4.5,-5) {Linear\\Datalog};
        \node (Leq2) at (9,-5) {\parbox{1.75cm}{\centering Horn-3SAT $\TLinP$}};
        \node[align=left] (HM2) at (13.5,-5) {Kearnes- \\Kiss};
        
        \node (L3) at (0,-8) {$\in \operatorname{P}$};
        \node (Leq3) at (9,-8) {\centering $\digraph K_3$};
        \node[align=left] (HM3) at (13.5,-8) {Kearnes-\\Markovi\'{c}-\\McKenzie}; 
        
        \draw (-1.7,-0.7) -- (16,-0.7);
        
        \path
        (HM1) edge[<->] node[above] {HMcK} (Leq1)
        (LSD1) edge[->] node[below] {E} (Leq1)
        (LSD1) edge[->] node[above] {LT} (L1)
        
        (HM2) edge[<->] node[above] {HMcK} (Leq2)
        (LSD2) edge[->] node[below] {AC} (Leq2)
        (LSD2) edge[->] node[above] {D} (L2)
        
        (HM3) edge[<->] node[above] {KMM} (Leq3)
        
        (L1) edge[->] (L2)
        (L2) edge[->] (L3)
        (LSD1) edge[->] (LSD2)
        (Leq1) edge[->] (Leq2)
        (Leq2) edge[->] (Leq3)
        (HM1) edge[->] (HM2)
        (HM2) edge[->] (HM3)
        
        (L3) edge[->,bend right=10,dashed] node[below] {$\operatorname{P}\neq\operatorname{NP}$} (Leq3)        
        (Leq3) edge[->,bend right=10] node[above] {BZ} (L3)

        (L2) edge[->,bend right=18,dashed] node[below] {$\operatorname{NL}\neq\operatorname{P}$, $\operatorname{NL}\neq\operatorname{Mod}_p\text{L}$} (Leq2)
        (L1) edge[->,bend right=18,dashed] node[below] {$\operatorname{L}\neq\operatorname{NL}$, $\operatorname{L}\neq\operatorname{Mod}_p\text{L}$} (Leq1)
        ;
        \path[dashed]
        
        (Leq1) edge[->,bend right=15,in=190] node[above] {Conj.~\ref{conj:L}} (LSD1)
        (Leq2) edge[->,bend right=15] node[above] {Conj.~\ref{conj:NL}} (LSD2)
        
        ;
    \end{tikzpicture}
    \caption{An overview of (computational and descriptive) complexity classes that are relevant for finite-domain CSPs, of important pp-constructions, and of the respective polymorphism conditions.
    LT stands for Larose and Tesson~\cite{LaroseTesson}, E stands for Egri \cite{EgriLT08}, HMcK stands for Hobby-McKenzie~\cite{HobbyMcKenzie} (Chapter 9), D stands for Dalmau~\cite{LinearDatalog}, 
    AC stands for Afrati and Cosmadakis~\cite{AfratiCosmadakis}, BZ stands for Bulatov~\cite{BulatovFVConjecture} and Zhuk~\cite{ZhukFVConjecture}, and KMM stands for Kearnes, Markovi\'c, and McKenzie~\cite{KearnesMarkovicMcKenzie}. 
    }
    \label{fig:landscape}
\end{figure}

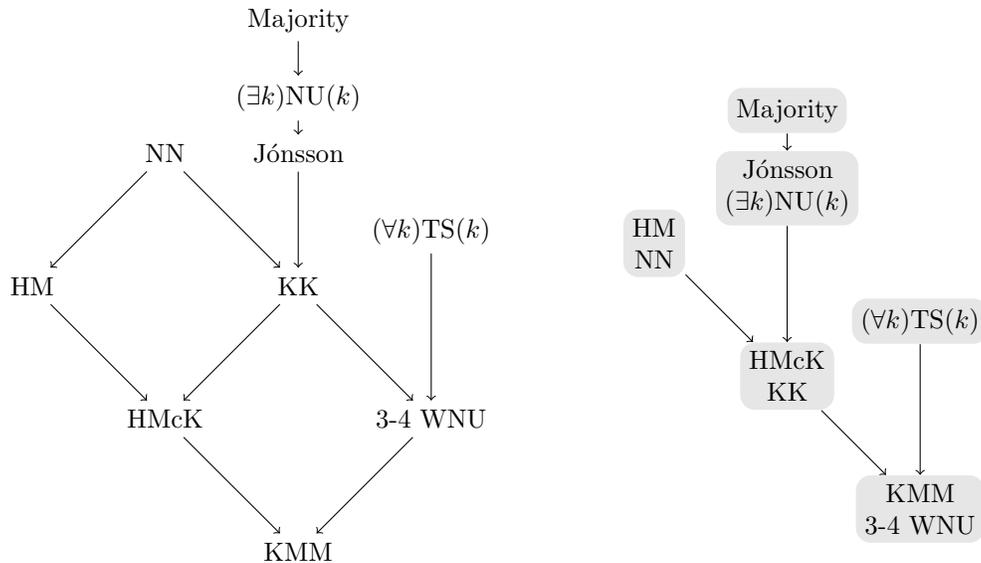
\begin{figure}
    \centering
    \begin{tikzpicture}[align=center, node distance=2.5cm]
        \node (34WNU) at (0,0) {3-4 WNU};
        \node[below left of = 34WNU] (KMM) {KMM};
        \node[above left of = 34WNU] (KK) {KK};
        \node[above left of = KMM] (HMcK) {HMcK};
        \node[above left of = HMcK] (HM) {HM};
        \node[above left of = KK] (NN) {NN};
        \node[above right of = NN] (Maj) {Majority};
        \node at (NN -| KK) (Jon) {J\'onsson};
        \node[above of = KK] (NU) {$(\exists k)\mathrm{NU}(k)$};
        \node[above of = 34WNU] (TS) {$(\forall k)\mathrm{TS}(k)$};
        
        \path[->]
        (Maj) edge (NU)
        (NU) edge (Jon)
        (Jon) edge (KK)
        (NN) edge (KK)
        (KK) edge (34WNU)
        (NN) edge (HM)
        (34WNU) edge (KMM)
        (HMcK) edge (KMM)
        (HM) edge (HMcK)
        (KK) edge (HMcK)
        (TS) edge (34WNU)
        ;
        
        \tikzset{class/.style={fill=black!10, minimum size=0.6cm, rounded corners=0.2cm}}
        
        \node[class] (34WNU2) at (6.5,-1.2) {KMM \\ 3-4 WNU};
        \node[class, above left of = 34WNU2] (KK2) {HMcK \\ KK};
        \node[class, above left of = KK2] (NN2) {HM \\ NN};
        \node[class, above right of = NN2] (Maj2) {Majority};
        \node[class, above of = KK2] (NU2) {J\'onsson \\ $(\exists k)\mathrm{NU}(k)$};
        \node[class, above of = 34WNU2] (TS2) {$(\forall k)\mathrm{TS}(k)$};
        
        \path[->]
        (Maj2) edge (NU2)
        (NU2) edge (KK2)
        (NN2) edge (KK2)
        (KK2) edge (34WNU2)
        (TS2) edge (34WNU2)
        ;
    \end{tikzpicture}
    
    \caption{
    Taxonomy of polymorphism conditions of structures with a finite domain that are relevant for the computational complexity of CSPs, ordered by strength; the arrows point from stronger conditions to weaker ones. The right picture shows the situation for trees assuming Conjecture~\ref{conj:bulin} (and P $\neq$ NP). 
    }
    
    \label{fig:tax}
\end{figure}

In this section we focus on certain concrete 
linear conditions that are relevant for studying the membership of CSPs in the most prominent complexity classes in the subsequent sections. 
An overview of the classes and the respective linear conditions 
is given in Figure~\ref{fig:landscape}. Solid arrows indicate implications, dotted arrows indicate conjectures. 
Figure~\ref{fig:tax} shows
the relationships between relevant linear polymorphism conditions that are defined throughout the section. The left side shows the general case and the right side shows the case for trees assuming Conjecture \ref{conj:bulin} (and $\operatorname{P}\neq\operatorname{NP}$). The implications are either immediate or from the literature~\cite{HobbyMcKenzie} (Chapter 9), ~\cite{BartoKozikWillard,KearnesMarkovicMcKenzie,Barto-cd}.

\subsection{Containment in P}
\label{sect:Siggers}
As discussed above, the characterization of the algebraic condition for tractability which is the most suitable 
for testing with a computer 
consists 
of a pair of ternary operations~\cite{KearnesMarkovicMcKenzie}. 

\begin{definition}
A pair of ternary operations $p,q \colon D^3 \to D$ is called \emph{Kearnes-Markovi\'{c}-McKenzie} if it satisfies the height-one condition
\begin{align*}
p(x,y,y) &\approx q(y,x,x)\approx q(x,x,y)\\
p(x,y,x)&\approx q(x,y,x).
\end{align*}
\end{definition}

Using this characterization, the CSP dichotomy can be stated as follows.

\begin{theorem}[\cite{ZhukFVConjecture,BulatovFVConjecture,KearnesMarkovicMcKenzie}]
A finite digraph $\HH$ has Kearnes-Markovi\'{c}-McKenzie polymorphisms if and only if there is no pp-construction of $\digraph K_3$ from $\HH$. In this case, $\Csp(\HH)$ is in P.
\end{theorem}

This characterization is optimal in the following sense: every height-one condition equivalent to Kearnes-Markovi\'{c}-McKenzie polymorphisms involves either an operation of arity at least 4 or at least two operations of arity 3 \cite{KearnesMarkovicMcKenzie}. 
However, there are several height-one conditions that imply 
the existence of Kearnes-Markovi\'{c}-McKenzie polymorphisms
and that are easier to test. In particular, we use the following. 

\begin{definition}\label{def:wnu}
An operation $f$ of arity $k$, for $k \geq 2$, is called a \emph{$k$-ary weak near-unanimity operation} (short, \emph{$k$-wnu}) if it satisfies the following  height-one condition 
$$f(y,x,\dots,x) \approx f(x,y,x,\dots,x) \approx \cdots \approx f(x,\dots,x,y).$$
A binary operation $f$ is called \emph{symmetric}
if it is a $2$-wnu, i.e., if it satisfies
$f(x,y) \approx f(y,x).$
\end{definition}

It is known that the existence of a $k$-wnu implies the existence of Kearnes-Markovi\'{c}-McKenzie polymorphisms
\cite{Siggers, KearnesMarkovicMcKenzie}, and that the existence of Kearnes-Markovi\'{c}-McKenzie polymorphisms 
implies the existence of a $k$-wnu for some $k \geq 2$~\cite{wnuf}.
Hence, in particular, if a finite digraph $\HH$ has a binary symmetric polymorphism then
$\Csp(\HH)$ can be solved in polynomial time. Our results will show that the converse is false even if $\HH$ is a tree, see Section~\ref{sect:NAC}. 

For both Kearnes-Markovi\'{c}-McKenzie and $k$-wnu, it is enough to test for level-wise satisfiability as discussed in Section \ref{sect:level-wise}. We prove the following more general claim.

\begin{lemma}\label{lemma:level-wise-two-variable}
Let $\Sigma$ be a height-one condition in two variables such that both the variables appear on each side in every identity from $\Sigma$. Then a balanced digraph level-wise satisfies $\Sigma$ if and only if it satisfies $\Sigma$.
\end{lemma}
\begin{proof}
Fix some polymorphisms $f,\dots$ that level-wise satisfy $\Sigma$, and define polymorphisms $f',\dots$ in the following way (say $f$ is $k$-ary): If $\level{x_1}=\level{x_2}=\dots=\level{x_k}$, define  $f'(x_1,x_2,\dots,x_k)=f(x_1,x_2,\dots,x_k)$. Else, let $\ell=\min\{\level{x_i}\mid 1\leq i\leq k\}$ and define $f'(x_1,x_2,\dots,x_k)=x_j$ where $j\in\{1,2,\dots,k\}$ is the smallest index such that $\level{x_j}=\ell$. 
To verify that the $f'$s are polymorphisms, note that if $(x_i,y_i)$ is an edge for $i\in\{1,2,\dots,k\}$, then $f'(x_1,x_2,\dots,x_k)$ and $f'(y_1,y_2,\dots,y_k)$ fall under the same case of the definition.
 If it is the second case, $x_j$ lies on the smallest level out of $\{\level{x_i}\mid 1\leq i\leq k\}$ if and only if $y_j$ lies on the smallest level out of $\{\level{y_i}\mid 1\leq i\leq k\}$. Hence, the selected coordinate $j$ is the same. To see that every identity is satisfied, note that the only interesting case is when $\level{x}\neq\level{y}$, and $f'$ then chooses the variable on the lower level. The other implication is trivial.
\end{proof}

\subsection{Containment in Datalog}
We have already mentioned in the introduction that containment in Datalog has numerous equivalent characterizations. In this section, we formally state one of these characterizations in terms of pp-constructibility and one in terms of height-one conditions. 
The structure $\TLinP$ has the domain
$D = \{0,\dots,p-1\}$ where $p$ is some prime, the relation $\{(x,y,z) \mid x+y+z\equiv 0 \pmod p\}$, 
and the relation $\{x\}$ for every $x \in D$. 
It is well-known that 
$\Csp(\TLinP)$ is not in Datalog~\cite{FederVardi}. 

\begin{definition}\label{def:wnu34}
A \emph{3-4 weak near-unanimity pair} (short, \emph{3-4 WNU}) is a pair of operations $f,g$ such that $f$ is a 3-wnu, $g$ is a 4-wnu, and they additionally satisfy the identity
$$ f(x,x,y) \approx g(x,x,x,y).$$
\end{definition}

\begin{theorem}[\cite{BoundedWidthJournal,Maltsev-Cond}]
Let $\HH$ be a finite digraph. Then the following are equivalent. 
\begin{itemize}
\item $\HH$ can be solved by Datalog.
\item there is no pp-construction of $\TLinP$ in $\HH$, for any prime $p$. 
\item $\HH$ has a 3-4 weak near-unanimity pair of polymorphisms.
\end{itemize}
\end{theorem}

Note that this shows that Conjecture~\ref{conj:bulin} implies (assuming P $\neq$ NP) that every tree with Kearnes-Markovi\'{c}-McKenzie polymorphisms has a 3-4WNU pair of polymorphisms and in particular, that it has a 3-wnu polymorphism, which is open as well. Also note that by Lemma \ref{lemma:level-wise-two-variable}, it is enough to test for  level-wise 3-4 WNU.

\subsection{Containment in NL}
\label{sect:NL}
In this section we present a strong sufficient condition 
for the containment of $\HH$ in NL.

\begin{definition}\label{def:j} 
For $n \geq 0$, a \emph{J\'onsson chain of length $n$ over $D$} is a sequence
of ternary operations  $j_1,j_2,\dots,j_{2n+1}$ on $D$ 
that satisfy
\begin{align*}
x & \approx j_1(x,x,y) \\
j_{2i-1}(x,y,y) & \approx j_{2i}(x,y,y) && \text{ for all } i \in \{1,\dots,n\} \\
j_i(x,y,x) & \approx x && \text{ for all } i \in \{1,\dots,2n+1\} \\
j_{2i}(x,x,y) & \approx j_{2i+1}(x,x,y) && \text{ for all } i  \in \{1,\dots,n\} \\
j_{2n+1}(x,y,y) & \approx y.
\end{align*}
The respective height one condition is abbreviated by $\J(n)$.
\end{definition}

Note that $\J(n)$ implies $\J(n+1)$ for every $n \geq 0$. Also note that for $n=0$ the operation $j_1$ must be a so-called \emph{majority} operation, which is the ternary case of a \emph{near-unanimity (NU)} operation, that is, an operation satisfying the identities 
$$x\approx f(y,x,\dots,x) \approx f(x,y,x,\dots,x) \approx \cdots \approx f(x,\dots,x,y).$$

The existence of a near-unanimity polymorphism characterizes \emph{bounded strict width} \cite{FederVardi}. More importantly for us, a near-unanimity polymorphism is sufficient to put $\HH$ in NL, using the following two results. Barto, Kozik, and Willard proved that finite structures with finite relational signature and a near-unanimity polymorphism have bounded pathwidth duality~\cite{BartoKozikWillard}. Dalmau proved that bounded pathwidth duality implies containment in NL~\cite{LinearDatalog}.

Barto~\cite{Barto-cd} moreover proved that if a finite structure with a finite relational signature
has polymorphisms that form a J\'onsson chain, then it also has a near-unanimity polymorphism (albeit its arity in the proof is doubly exponential in the size of the domain). In the other direction, it is well known that $\NU(n)$ implies $\J(n-2)$, syntactically. Therefore, we do not test for near-unanimities of arities higher than 3; it is more efficient to test for a J{\'o}nsson chain.

\begin{theorem}[\cite{Barto-cd,BartoKozikWillard,LinearDatalog}]
If a finite digraph $\HH$ satisfies $\J(n)$ for some $n \geq 1$, 
then $\Csp(\HH)$ is in linear Datalog, and hence in NL. 
\end{theorem}

Note that the existence of polymorphisms of $\HH$ that form a J\'onsson chain is only a sufficient condition for the containment of $\Csp(\HH)$ in NL. An incomparable sufficient condition for the containment of $\Csp(\HH)$ in NL was identified in~\cite{CarvalhoDalmauKrokhin}. The condition presented there also has a characterization via height-one identities, but the arities of the operations are prohibitively large so that we did not implement this test for trees.

The conjectured characterization of containment in NL, the Kearnes-Kiss chain from Conjecture \ref{conj:NL}, is defined in Section \ref{sect:PorModpL-hardness} (Definition \ref{def:kk}).

\subsection{Containment in L}
\label{sect:L}
One of the strongest known sufficient conditions for containment in L is a conditional result of Kazda, which involves the following linear condition.

\begin{definition}\label{def:hm} 
For $n \geq 1$, a \emph{Hagemann-Mitschke chain of length $n$ over $D$} is a sequence
of ternary operations  $p_1,\dots,p_{n}$ on $D$ 
that satisfy
\begin{align*}
x & \approx p_1(x,y,y) \\
p_{i}(x,x,y) & \approx p_{i+1}(x,y,y) && \text{ for all } i \in \{1,\dots,n-1\} \\
p_n(x,x,y) & \approx y.
\end{align*}
The respective height one condition is abbreviated by $\HM(n)$.
\end{definition}

Note that $\HM(n)$ implies $\HM(n+1)$ for every $n \geq 1$.
For $n=1$ the operation $p_1$ is known as a Maltsev operation. Kazda~\cite{Kazda-n-permute} proved the following conditional result.

\begin{theorem}[\cite{Kazda-n-permute}]
If a finite digraph $\HH$ can be solved by linear Datalog, and 
$\HH$ satisfies $\HM(n)$ for some $n \geq 1$, 
then $\HH$ can also be solved by linear symmetric Datalog (and hence is in L).
\end{theorem}

The conjectured characterization of containment in L, the Noname chain from Conjecture \ref{conj:L}, is defined in Section \ref{sect:NLorModpL-hardness} (Definition \ref{def:nn}).

\subsection{Solvability by Arc Consistency}
\label{sect:AC-Pol}
Solvability by Arc-Consistency
(and tree duality) can be characterized in terms of height one conditions as well. 

\begin{definition}
An operation $s_n \colon D^n \to D$ is called \emph{totally symmetric} if for all variables
$x_1,\dots,x_n$ and $y_1,\dots,y_n$ such that $\{x_1,\dots,x_n\} = \{y_1,\dots,y_n\}$ the operation $s_n$ satisfies 
$$s_n(x_1,\dots,x_n) \approx s_n(y_1,\dots,y_n).$$
The respective height one condition is abbreviated by $\TS(n)$. 
\end{definition}

The digraph $\HH$ can be solved by arc consistency if and only if $\HH$ has totally symmetric polymorphisms of all arities~\cite{FederVardi,DalmauPearson}. Note that $\TS(4)$ implies 3-4 WNU. Also note that a finite digraph $\HH$ satisfies $\TS(n)$ for all $n>0$ if and only if it satisfies $\TS(2 |E(\HH)|)$ (see the proof given in~\cite{DalmauPearson}). The arity $2 |E(\HH)|$ is still fairly large; therefore it is particularly useful that the level-wise test is sufficient.

\begin{lemma}[see {\cite[proof of Lemma 4.1]{SpecialTriads}}]\label{lemma:level-wise-ts}
For any balanced digraph $\HH$ and $n>0$, $\HH$ level-wise satisfies $\TS(n)$  if and only if $\HH$ satisfies $\TS(n)$.
\end{lemma}
\begin{proof}
Let $s_n$ be an $n$-ary polymorphism of $\HH$ that level-wise satisfies the condition $\TS(n)$. 
We can construct an $n$-ary totally symmetric polymorphism $s'_n$ of $\HH$ by applying $s_n$ to the set of vertices on the smallest level. That is, for an input tuple $(x_1,\dots,x_n)$ let $\ell=\min\{\level{x_i}\mid 1\leq i\leq n\}$, $\{x_i\mid\level{x_i}=\ell\}=\{x_{i_1},\dots,x_{i_k}\}$, and set 
\[
s'_n(x_1,\dots,x_n)=s_n(x_{i_1},\dots,x_{i_k},\underbrace{x_{i_k},\dots,x_{i_k}}_{(n-k)\text{ times}}).
\]
Clearly, the definition of $s'_n(x_1,\dots,x_n)$ depends only on the set $\{x_1,\dots,x_n\}$. To see that $s'_n$ is a polymorphism, note that similarly as in Lemma \ref{lemma:level-wise-two-variable} if $(x_i,y_i)$ is an edge for $i\in\{1,2,\dots,n\}$, then $x_{i_j}$ lies on the smallest level if and only if $y_{i_j}$ does. The rest follows from the fact that $s_n$ is totally symmetric on each level. The other implication is again trivial.
\end{proof}

\subsection{P-hardness}
\label{sect:P-hardness}
The structure $\HornSAT$
has the domain $\{0,1\}$ and a ternary relation $\{0,1\}^3 \setminus \{(1,1,0)\}$,
and the two unary relations $\{0\}$ and $\{1\})$.
It is well-known that $\Csp(\HornSAT)$ is P-complete, i.e., complete for the complexity class P under deterministic logspace reductions.

\begin{definition}
A \emph{Hobby-McKenzie chain} consists of ternary operations $d_0,\dots,d_n,p,e_0,\dots,e_n$
such that 
\begin{align*}
    d_0(x,y,z) & \approx x \\
    d_i(x,y,y) & \approx d_{i+1}(x,y,y) & \text{for even } i < n \\
    d_i(x,x,y) & \approx d_{i+1}(x,x,y) & \text{for odd } i < n \\
    d_i(x,y,x) & \approx d_{i+1}(x,y,x) & \text{for odd } i < n \\
    d_n(x,y,y) & \approx p(x,y,y) \\
    p(x,x,y) & \approx e_0(x,x,y) \\
    e_i(x,y,y) & \approx e_{i+1}(x,y,y) & \text{for even } i < n \\
    e_i(x,x,y) & \approx e_{i+1}(x,x,y) & \text{for odd } i < n \\
    e_i(x,y,x) & \approx e_{i+1}(x,y,x) & \text{for even } i < n \\
    e_n(x,y,z) & \approx z.
\end{align*}
\end{definition}
The respective height one condition is abbreviated by $\HMcK(n)$.

\begin{theorem}[consequence of Theorem 9.8 in~\cite{HobbyMcKenzie}]
A finite structure does not satisfy $\HMcK(n)$ for some $n \geq 1$ if and only if it can pp-construct $\HornSAT$.
\end{theorem}

The theorem implies that if a finite digraph $\HH$ does not satisfy $\HMcK(n)$ for some $n \geq 1$, then $\HH$ is P-hard (see Section~\ref{sect:pp}).

\subsection{P-hardness or Mod$_p$L-hardness} \label{sect:PorModpL-hardness}
It is widely believed that NL is a proper subclass of P. 
Another complexity class which is believed to be a proper subclass of P is the class 
Mod$_p$L, for some prime $p$: this is defined to be the class of problems such that there exists a non-deterministic logspace machine $M$ such that an instance is in the class if and only if the number of accepting paths of $M$ on the instance is divisible by $p$; see~\cite{LaroseTesson}. It is well known that $\Csp(\TLinP)$ is Mod$_p$L-complete. 
If NL would contain Mod$_p$L then this would be a considerable breakthrough in complexity theory.

\begin{definition}[from Theorem 9.11 in~\cite{HobbyMcKenzie}] \label{def:kk}
Let $D$ be a set. A \emph{Kearnes-Kiss chain of length
$n \geq 2$ over $D$} is a sequence
of ternary operations  $d_0,d_1,\dots,d_n$ on $D$ such that 
\begin{align*}
d_0(x,y,z) & \approx x \\
d_i(x,y,y) & \approx d_{i+1}(x,y,y) && \text{ for even } i \in \{0,1,\dots,n-1\} 
\\ 
d_i(x,y,x) & \approx d_{i+1}(x,y,x) 
&& \text{ for even } i \in \{0,1,\dots,n-1\} \\
d_i(x,x,y) & \approx d_{i+1}(x,x,y) && \text{ for odd } i  \in \{1,\dots,n-1\} \\
d_n(x,y,z) & \approx z .
\end{align*}
The respective height one condition is abbreviated by $\KK(n)$.
\end{definition}
Note that $\KK(n)$ implies $\KK(n+1)$ for every $n \geq 0$.  
Also note that the existence of a J{\'o}nsson chain implies the existence of a Kearnes-Kiss chain \cite{HobbyMcKenzie}; namely $
\J(n)$ trivially implies $\KK(2n+4)$.

\begin{theorem}[see~\cite{BartoKozikWillardPolymorphisms} and~\cite{HobbyMcKenzie}]
A finite structure does not satisfy $\KK(n)$ for some $n \geq 2$ if and only if it can pp-construct $\HornSAT$ or $\TLinP$ for some prime $p$. 
\end{theorem}

 We have already mentioned that if a finite digraph $\HH$ pp-constructs $\HornSAT$ then it is P-hard, and hence $\Csp(\HH)$ is not in NL unless $\operatorname{NL}=\operatorname{P}$. Similarly, 
 if $\HH$ pp-constructs $\TLinP$ then it is Mod$_p$L-hard, and in this case
it is not in NL unless NL contains Mod$_p$L.
If the conjecture that
`easy trees cannot count' (Conjecture~\ref{conj:bulin}) is true, then the existence a Kearnes-Kiss chain is equivalent to the existence of a Hobby-McKenzie chain for trees (assuming P $\neq$ NP).

\subsection{NL-hardness} \label{sect:NL-hardness}
The structure $\stcon$ has the domain $\{0,1\}$, the  binary relation $\{(0,0),(0,1),(1,1)\}$ and the unary relations $t = \{0\}$ and $s = \{1\}$. Note that an instance of the CSP of $\stcon$ is unsatisfiable if and only if there exists a directed path from $s$ to $t$ in the digraph defined by the binary relation. 
It is well-known that $\Csp(\stcon)$ is complete for the complexity class NL. 

\begin{theorem}[\cite{LaroseTesson,HobbyMcKenzie}]
If a finite digraph $\HH$ does not satisfy $\HM(n)$ for some $n \geq 1$, then it can pp-construct the structure $\stcon$, and $\HH$ is NL-hard. 
\end{theorem}

\subsection{NL-hardness or Mod$_p$L-hardness} \label{sect:NLorModpL-hardness}
We now present a polymorphism condition that characterizes the finite structures that can pp-construct 
$\stcon$ \emph{or} $\TLinP$ for some prime $p$.

\begin{definition}[\cite{HobbyMcKenzie}, Theorem 9.15]\label{def:nn} 
For $n \geq 0$, a \emph{Noname chain of length $n$ over $D$} is a sequence
of operations  $f_0,f_1,\dots,f_n$ of arity four on $D$ such that 
\begin{align*}
f_0(x,y,y,z) & \approx x \\
f_i(x,x,y,x) & \approx f_{i+1}(x,y,y,x) && \text{for all } i \in \{0,\dots,n-1\} \\
f_i(x,x,y,y) & \approx f_{i+1}(x,y,y,y)
&& \text{for all } i \in \{0,\dots,n-1\} \\
f_n(x,x,y,z) & \approx z. 
\end{align*}
The respective height one condition is abbreviated by $\NN(n)$.
\end{definition}
Note that $\NN(n)$ implies $\NN(n+1)$ for every $n \geq 0$.

\begin{theorem}[\cite{HobbyMcKenzie, BartoKozikWillardPolymorphisms}]
A finite structure does not satisfy $\NN(n)$ for some $n \geq 1$ if and only if it can pp-construct the structure $\stcon$ or the structure $\TLinP$ for some prime $p$. 
\end{theorem}

It follows that if a finite digraph $\HH$ does not satisfy $\NN(n)$ for any $n \geq 1$, then $\Csp(\HH)$ is NL-hard or Mod$_p$L-hard.
Hence, $\HH$ is in this case not in L, unless $\operatorname{L}=\operatorname{NL}$ or L = Mod$_p$L.
Note that Conjecture~\ref{conj:bulin}\blue together with Conjecture~\ref{conj:L} implies that $\NN(n)$ for some $n$ and $\HM(n)$ for some $n$ are equivalent for trees (assuming $\operatorname{L}\neq\operatorname{NL}$).

\section{Experimental Results}
\label{sect:exp}

We implemented the AC-3 algorithm for establishing arc-consistency and used its adaptation, known as the MAC-3 algorithm, 
for maintaining arc-consistency during the backtracking procedure described in Section~\ref{subsect:testing-linear}. 
The lists and related operations were implemented by doubly-linked lists. The code is written in 
Rust\footnote{https://gitlab.com/WhatDothLife/tripolys} and the experiments were run on a Intel(R) Xeon(R) CPU E5-2680 v3 (12 cores) @ 2.50GHz with Linux.
We also used another implementation written in Python. All tests for chains of polymorphisms and for totally symmetric polymorphisms where done using this implementation on a AMD Ryzen 5 4500U (with 8 cores) @ 2.38 GHz with Windows. 
An efficient implementation was essential to obtain our results.\footnote{For comparison, we modeled some polymorphism tests in the constraint modeling language MiniZinc \cite{minizinc} and used the CP solver Gecode \cite{gecode}. However, this turned out to not be competitive.}

Table~\ref{table:otrees} shows the number of unlabeled trees with $n$ vertices and the number of those that are cores.
The table suggests that the
fraction of trees that are cores quickly goes to 0. The next columns contain the number of unlabeled rooted cores with $n$ vertices, the number of core checks, 
and the mean cpu time per core check on a tree with $n$ vertices. The final column in the table shows the computation time needed to generate all the unlabeled 
core trees with $n$ vertices with Algorithm \ref{alg:genTrees}.

\begin{table}
\begin{center}
\begin{tabular}{rrrrrc S[table-format=3.5]}
\toprule
\thead{$n$} & \thead{trees} & \thead{cores} & \thead{rooted \\ cores} & \thead{core \\ checks} & \thead{time per \\ core check (\si{\micro\second})} & {\thead{total time}}\\
\toprule
1   & 1            & 1       & 1          & 1       &  10      &        0.7\,\si{\milli\second}  \\
2   & 1            & 1       & 2          & 1       &  10      &        0.8\,\si{\milli\second}  \\
3   & 3            & 1       & 3          & 3       &  11      &        1.0\,\si{\milli\second}  \\
4   & 8            & 1       & 6          & 8       &  12      &        1.1\,\si{\milli\second}  \\
5   & 27           & 1       & 11         & 19      &  13      &        1.3\,\si{\milli\second}  \\
6   & 91           & 2       & 28         & 39      &  15       &        2.0\,\si{\milli\second}  \\
7   & 350          & 3       & 63         & 94      &  24       &        2.3\,\si{\milli\second}  \\
8   & 1376         & 7       & 170        & 198     &  25       &        5.6\,\si{\milli\second}  \\
9   & 5743         & 15      & 439        & 439     &  30       &       12.0\,\si{\milli\second}  \\
10  & 24635        & 36      & 1200       & 953     &  36       &       30.9\,\si{\milli\second}  \\
11  & 108968       & 85      & 3307       & 2180    &  40       &       69.1\,\si{\milli\second}  \\
12  & 492180       & 226     & 9380       & 5050    &  48       &      187.9\,\si{\milli\second}  \\
13  & 2266502      & 578     & 26731      & 12218   &  55       &      547.2\,\si{\milli\second}  \\
14  & 10598452     & 1569    & 77508      & 29785   &  66       &        1.6\,\si{\second}  \\
15  & 50235931     & 4243    & 226399     & 74902   &  71      &    5.0\,\si{\second}  \\
16  & 240872654    & 11848   & 668228     & 190632  &  84      &   15.9\,\si{\second}  \\
17  & 1166732814   & 33104   & 1984592    & 496373  &  98      &   52.5\,\si{\second}  \\
18  & 5702001435   & 94221   & 5937276    & 1308847 &  121      &    3.0\,\si{\minute}  \\
19  & 28088787314  & 269455  & 17856807   & 3512229 &  129      &    9.8\,\si{\minute}  \\
20  & 139354922608 & 779268  & 53996424   & 9538804 &  142        &   31.7\,\si{\minute}
\end{tabular}
\end{center}
\caption{The number of unlabeled trees with $n$ vertices together with results of Algorithm~\ref{alg:genTrees}.}
\label{table:otrees}
\end{table}

In this section we present the results of testing the discussed linear conditions on these trees to classify them with respect to their computational complexity. In some cases we manage to compute all the minimal trees in the respective complexity class; the corresponding results are presented in Section~\ref{sect:sht}. In Section~\ref{sect:opentrees} we present trees whose precise complexity status is open. 
While the numbers of trees given in the text are up to isomorphism, the corresponding figures make a further restriction based on the following fact.

\begin{remark}
\label{rem:reverse}
An operation is a polymorphism of $\HH$ if and only if it is a polymorphism of $\HH^R$.
\end{remark}

The remark justifies that 
our figures contain exactly one of the trees $\T$, $\T^R$.
It turns out that the trees in our figures that do not satisfy a certain height-one condition have a unique minimal subtree which has no \emph{idempotent} polymorphisms satisfying the respective condition (see Remark \ref{rem:idempotence}).
The vertices and edges drawn in gray do not belong to this minimal subgraph of $\T$.

\subsection{The Smallest Hard Trees}
\label{sect:sht}
In this section we present the smallest trees that are NP-hard and that are NL-hard, under standard assumptions from complexity theory. 
We also compute the smallest tree that cannot be solved by arc consistency, the smallest trees that cannot be solved by Datalog,
and the smallest trees that cannot be solved by linear symmetric Datalog; these results hold without any assumptions from complexity theory.

\subsubsection{The Smallest NP-Hard Trees}
Our algorithm found that all trees with at most 19 vertices have 
Kearnes-Markovi\'{c}-McKenzie polymorphisms 
and hence are tractable. 
It also found that there exist exactly 36 trees with 20 vertices that have no Kearnes-Markovi\'{c}-McKenzie polymorphisms
and hence are NP-hard. 
For such an NP-hard tree $\T$ with 20 vertices it takes our algorithm about 0.07 seconds to construct the indicator digraph $\indicator{\T}$ for the Kearnes-Markovi\'{c}-McKenzie polymorphisms and about 0.03 seconds to verify that $\indicator{\T}$ does not have a homomorphism to $\T$. When applying the level trick, it takes about 0.01 seconds to construct the indicator digraph and 0.03 seconds to verify that $\indicator{\T}$ does not have a homomorphism to $\T$.

The trees with 20 vertices that have no Kearnes-Markovi\'{c}-McKenzie polymorphisms 
are displayed in Figure~\ref{fig:hard-trees}.
Note that the smallest subtree without idempotent Kearnes-Markovi\'{c}-McKenzie polymorphism is the same for the trees A2--A4, A5--A12, and A13--A18.

\input{trees}

Moreover, there are 4 smallest triads with 22 vertices that have no Kearnes-Markovi\'{c}-McKenzie polymorphisms;
these are shown in Figure~\ref{fig:hardTriads}. 
All smaller triads have a binary symmetric polymorphism.

\def\scale{0.4}
\def\hdist{1cm}
\def\vdist{1cm}
\begin{figure}
\centering
\begin{tikzpicture}[scale=\scale]
\node[bullet,gray] (21) at (4,0) {};
\node[bullet] (4) at (0.5,1) {};
\node[bullet] (14) at (6,1) {};
\node[bullet] (20) at (4,1) {};
\node[bullet] (1) at (2,2) {};
\node[bullet] (5) at (0,2) {};
\node[bullet] (3) at (1,2) {};
\node[bullet] (9) at (5,2) {};
\node[bullet] (13) at (6,2) {};
\node[bullet] (19) at (4,2) {};
\node[bullet] (2) at (1,3) {};
\node[bullet] (6) at (0,3) {};
\node[bullet] (0) at (2,3) {};
\node[bullet] (10) at (5,3) {};
\node[bullet] (12) at (6,3) {};
\node[bullet] (16) at (3,3) {};
\node[bullet] (18) at (4,3) {};
\node[bullet] (7) at (0,4) {};
\node[bullet] (11) at (5.5,4) {};
\node[bullet] (15) at (2.5,4) {};
\node[bullet] (17) at (3.5,4) {};
\node[bullet,gray] (8) at (0,5) {};
\path[->,>=stealth']
(1) edge (2)
(1) edge (0)
(4) edge (5)
(4) edge (3)
(5) edge (6)
(6) edge (7)
(7) edge[gray] (8)
(0) edge (15)
(3) edge (2)
(9) edge (0)
(9) edge (10)
(10) edge (11)
(12) edge (11)
(13) edge (12)
(14) edge (13)
(16) edge (15)
(16) edge (17)
(18) edge (17)
(19) edge (18)
(20) edge (19)
(21) edge[gray] (20)
;
\node at (3,-1) {Triad 1};
\end{tikzpicture}
\hspace{\hdist}
\begin{tikzpicture}[scale=\scale]
\node[bullet,gray] (21) at (6,0) {};
\node[bullet] (4) at (0.5,1) {};
\node[bullet] (14) at (6,1) {};
\node[bullet] (20) at (4,1) {};
\node[bullet] (1) at (2,2) {};
\node[bullet] (5) at (0,2) {};
\node[bullet] (3) at (1,2) {};
\node[bullet] (9) at (5,2) {};
\node[bullet] (13) at (6,2) {};
\node[bullet] (19) at (4,2) {};
\node[bullet] (2) at (1,3) {};
\node[bullet] (6) at (0,3) {};
\node[bullet] (0) at (2,3) {};
\node[bullet] (10) at (5,3) {};
\node[bullet] (12) at (6,3) {};
\node[bullet] (16) at (3,3) {};
\node[bullet] (18) at (4,3) {};
\node[bullet] (7) at (0,4) {};
\node[bullet] (11) at (5.5,4) {};
\node[bullet] (15) at (2.5,4) {};
\node[bullet] (17) at (3.5,4) {};
\node[bullet,gray] (8) at (0,5) {};
\path[->,>=stealth']
(1) edge (2)
(1) edge (0)
(4) edge (5)
(4) edge (3)
(5) edge (6)
(6) edge (7)
(7) edge[gray] (8)
(0) edge (15)
(3) edge (2)
(9) edge (0)
(9) edge (10)
(10) edge (11)
(12) edge (11)
(13) edge (12)
(14) edge (13)
(16) edge (15)
(16) edge (17)
(18) edge (17)
(19) edge (18)
(20) edge (19)
(21) edge[gray] (14)
;
\node at (3,-1) {Triad 2};
\end{tikzpicture}
\caption{The smallest NP-hard triads (up to edge reversal and assuming P $\neq$ NP).}
\label{fig:hardTriads}
\end{figure}

\subsubsection{The Smallest NL-hard Trees}
There are 8 trees with 12 vertices that are NL-hard. 
Two of them are isomorphic to their reverse, so we only 
 display 5 trees in Figure~\ref{fig:NLHardTrees}, called B1, B2, B3, B4, and B5. The proof that they are NL-hard can be found below.
All other trees with at most 12 vertices satisfy  $\HM(8)$.

\def\scale{0.4}
\def\hdist{1cm}
\def\vdist{1cm}
\begin{figure}
\centering\begin{tikzpicture}[scale=\scale]
\node[bullet] (11) at (0,0) {};
\node[bullet] (6) at (2,0) {};
\node[bullet] (7) at (1,-1) {};
\node[bullet] (5) at (2,-1) {};
\node[bullet] (10) at (0,-1) {};
\node[bullet] (8) at (0,-2) {};
\node[bullet] (0) at (1,-2) {};
\node[bullet] (4) at (2,-2) {};
\node[bullet] (3) at (1,-3) {};
\node[bullet] (9) at (0,-3) {};
\node[bullet,gray] (1) at (2,-3) {};
\node[bullet,gray] (2) at (2,-4) {};
\path[<-,>=stealth']
(7) edge (8)
(7) edge (0)
(8) edge (9)
(0) edge (3)
(0) edge[gray] (1)
(5) edge (4)
(4) edge (3)
(1) edge[gray] (2)
(11) edge (10)
(10) edge (8)
(6) edge (5)
;
\node at (1,-5) {Tree B1};
\end{tikzpicture}
\hspace{\hdist}
\begin{tikzpicture}[scale=\scale]
\node[bullet,gray] (9) at (0,0) {};
\node[bullet] (0) at (1,1) {};
\node[bullet] (8) at (0,1) {};
\node[bullet] (7) at (0,2) {};
\node[bullet] (5) at (1,2) {};
\node[bullet] (1) at (2,2) {};
\node[bullet] (2) at (2,3) {};
\node[bullet] (10) at (0,3) {};
\node[bullet] (6) at (1,3) {};
\node[bullet] (3) at (2,4) {};
\node[bullet] (11) at (0,4) {};
\node[bullet,gray] (4) at (2,5) {};
\path[->,>=stealth']
(2) edge (3)
(3) edge[gray] (4)
(7) edge (10)
(7) edge (6)
(10) edge (11)
(0) edge (5)
(0) edge (1)
(5) edge (6)
(1) edge (2)
(9) edge[gray] (8)
(8) edge (7)
;
\node at (1,-1) {Tree B2};
\end{tikzpicture}
\hspace{\hdist}
\begin{tikzpicture}[scale=\scale]
\node[bullet,gray] (5) at (0,0) {};
\node[bullet] (4) at (0,-1) {};
\node[bullet] (11) at (2,-1) {};
\node[bullet] (0) at (1,-2) {};
\node[bullet] (3) at (0,-2) {};
\node[bullet] (10) at (2,-2) {};
\node[bullet] (6) at (1,-3) {};
\node[bullet] (9) at (2,-3) {};
\node[bullet] (1) at (0,-3) {};
\node[bullet] (7) at (1,-4) {};
\node[bullet] (2) at (0,-4) {};
\node[bullet,gray] (8) at (1,-5) {};
\path[<-,>=stealth']
(0) edge (6)
(0) edge (1)
(6) edge (7)
(4) edge (3)
(3) edge (1)
(7) edge[gray] (8)
(9) edge (7)
(1) edge (2)
(11) edge (10)
(10) edge (9)
(5) edge[gray] (4)
;
\node at (1,-6) {Tree B3};
\end{tikzpicture}
\hspace{\hdist}
\begin{tikzpicture}[scale=\scale]
\node[bullet,gray] (3) at (0,0) {};
\node[bullet] (2) at (0,1) {};
\node[bullet] (9) at (2,1) {};
\node[bullet] (6) at (1,2) {};
\node[bullet] (1) at (0,2) {};
\node[bullet] (8) at (2,2) {};
\node[bullet] (7) at (2,3) {};
\node[bullet] (4) at (0,3) {};
\node[bullet] (0) at (1,3) {};
\node[bullet] (10) at (2,4) {};
\node[bullet] (5) at (0,4) {};
\node[bullet,gray] (11) at (2,5) {};
\path[->,>=stealth']
(6) edge (7)
(6) edge (0)
(7) edge (10)
(1) edge (4)
(1) edge (0)
(4) edge (5)
(8) edge (7)
(2) edge (1)
(10) edge[gray] (11)
(9) edge (8)
(3) edge[gray] (2)
;
\node at (1,-1) {Tree B4};
\end{tikzpicture}
\hspace{\hdist}
\begin{tikzpicture}[scale=\scale]
\node[bullet,gray] (11) at (0,0) {};
\node[bullet] (10) at (0,-1) {};
\node[bullet] (2) at (2,-1) {};
\node[bullet] (1) at (2,-2) {};
\node[bullet] (9) at (0,-2) {};
\node[bullet] (6) at (1,-2) {};
\node[bullet] (3) at (2,-3) {};
\node[bullet] (0) at (1,-3) {};
\node[bullet] (7) at (0,-3) {};
\node[bullet] (8) at (0,-4) {};
\node[bullet] (4) at (2,-4) {};
\node[bullet,gray] (5) at (2,-5) {};
\path[<-,>=stealth']
(1) edge (3)
(1) edge (0)
(3) edge (4)
(7) edge (8)
(9) edge (7)
(11) edge[gray] (10)
(10) edge (9)
(4) edge[gray] (5)
(6) edge (0)
(6) edge (7)
(2) edge (1)
;
\node at (1,-6) {Tree B5};
\end{tikzpicture}
\caption{The smallest NL-hard trees (up to edge reversal and assuming L $\neq$ NL).}
\label{fig:NLHardTrees}
\end{figure}

Since the trees B1-B5 have a majority they are in NL. 
To prove that B1-B5 are NL-hard, we show that they can pp-construct st-Con. Hence we need to construct the three relations $\{0\}, \{1\},$ and $\{(0,0),(0,1),(1,1)\}$. First 
note that in a core tree $\T$ any singleton set is pp-definable from $\T$, since  $\End(\T)=\{\operatorname{id}_{T}\}$ by Theorem~\ref{thm:ac}. The following two graphs represent two pp-formulas $\phi_1(x,y)$ and $\phi_2(x,y)$. The filled vertices stand for existentially quantified variables.

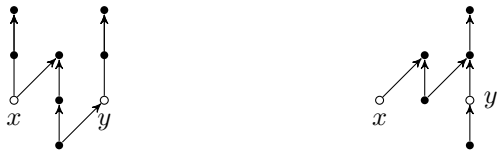
\begin{center}
\begin{tikzpicture}[scale=0.6]
    \node[bullethole,draw,label=below:$x$] (0) at (0,0) {};
    \node[bullet] (1) at (1,1) {};
    \node[bullet] (2) at (1,0) {};
    \node[bullet] (4) at  (1,-1) {};
    \node[bullethole,label=below:$y$] (6) at (2,0) {};
    \node[bullet] (5) at (0,2) {};
    \node[bullet] (9) at (0,1) {};
    \node[bullet] (7) at (2,2) {};
    \node[bullet] (8) at (2,1) {};
    
    \path[->,>=stealth']
        (0) edge (1)
        (2) edge (1)
        (4) edge (2)
        (4) edge (6)
        (9) edge (5)
        (0) edge (5)
        (8) edge (7)
        (6) edge (7)
        ;
\end{tikzpicture}
\hspace{3cm}
\begin{tikzpicture}[scale=0.6]
    \node[bullethole,draw,label=below:$x$] (0) at (0,0) {};
    \node[bullet] (1) at (1,1) {};
    \node[bullet] (2) at (1,0) {};
    \node[bullet] (3) at (2,1) {};
    \node[bullet] (4) at  (2,2) {};
    \node[bullethole,label=right:$y$] (6) at (2,0) {};
    \node[bullet] (7) at (2,-1) {};
    \path[->,>=stealth']
        (0) edge (1)
        (2) edge (1)
        (2) edge (3)
        (6) edge (3)
        (3) edge (4)
        (7) edge (6)
        ;
\end{tikzpicture}
\end{center}

The trees B1, B2, B3 can pp-define a structure that is homomorphically equivalent to $\stcon$ using  $\phi_1(x,y)$ for $\edges{\stcon}$. The trees B4 and B5 can do the same using  $\phi_2(x,y)$ for $\edges{\stcon}$. Since $\stcon$ is NL-hard, B1-B5 are NL-hard as well.

So for 12 vertices, 8 out of 226 trees are NL-hard (assuming L $\neq$ NL). In Table~\ref{fig:HMtrees} we present how this distribution in core trees changes with increasing number of vertices. Every tree with at most 20 vertices falls into one of two cases:
\begin{itemize}
    \item it satisfies $\HM(16)$ and has a majority polymorphism, hence it is in L, or
    \item it has no $\HM(30)$.
\end{itemize}

We strongly suspect that in the latter case the trees have no $\HM(n)$ for any $n$, can pp-construct st-Con, and are NL-hard.

\pgfplotstableread{ 
Label	1	2	3	4	5	6	7	8	9	10	11	12	13	14	15	16	17	18	19	20	21	22	23	24	25	26	27	28	29
1	1.00000	0.00000	0.00000	0.00000	0.00000	0.00000	0.00000	0.00000	0.00000	0.00000	0.00000	0.00000	0.00000	0.00000	0.00000	0.00000	0.00000	0.00000	0.00000	0.00000	0.00000	0.00000	0.00000	0.00000	0.00000	0.00000	0.00000	0.00000	0.00000
2	1.00000	0.00000	0.00000	0.00000	0.00000	0.00000	0.00000	0.00000	0.00000	0.00000	0.00000	0.00000	0.00000	0.00000	0.00000	0.00000	0.00000	0.00000	0.00000	0.00000	0.00000	0.00000	0.00000	0.00000	0.00000	0.00000	0.00000	0.00000	0.00000
3	1.00000	0.00000	0.00000	0.00000	0.00000	0.00000	0.00000	0.00000	0.00000	0.00000	0.00000	0.00000	0.00000	0.00000	0.00000	0.00000	0.00000	0.00000	0.00000	0.00000	0.00000	0.00000	0.00000	0.00000	0.00000	0.00000	0.00000	0.00000	0.00000
4	1.00000	0.00000	0.00000	0.00000	0.00000	0.00000	0.00000	0.00000	0.00000	0.00000	0.00000	0.00000	0.00000	0.00000	0.00000	0.00000	0.00000	0.00000	0.00000	0.00000	0.00000	0.00000	0.00000	0.00000	0.00000	0.00000	0.00000	0.00000	0.00000
5	1.00000	0.00000	0.00000	0.00000	0.00000	0.00000	0.00000	0.00000	0.00000	0.00000	0.00000	0.00000	0.00000	0.00000	0.00000	0.00000	0.00000	0.00000	0.00000	0.00000	0.00000	0.00000	0.00000	0.00000	0.00000	0.00000	0.00000	0.00000	0.00000
6	0.50000	0.50000	0.00000	0.00000	0.00000	0.00000	0.00000	0.00000	0.00000	0.00000	0.00000	0.00000	0.00000	0.00000	0.00000	0.00000	0.00000	0.00000	0.00000	0.00000	0.00000	0.00000	0.00000	0.00000	0.00000	0.00000	0.00000	0.00000	0.00000
7	0.33333	0.66667	0.00000	0.00000	0.00000	0.00000	0.00000	0.00000	0.00000	0.00000	0.00000	0.00000	0.00000	0.00000	0.00000	0.00000	0.00000	0.00000	0.00000	0.00000	0.00000	0.00000	0.00000	0.00000	0.00000	0.00000	0.00000	0.00000	0.00000
8	0.14286	0.71429	0.00000	0.14286	0.00000	0.00000	0.00000	0.00000	0.00000	0.00000	0.00000	0.00000	0.00000	0.00000	0.00000	0.00000	0.00000	0.00000	0.00000	0.00000	0.00000	0.00000	0.00000	0.00000	0.00000	0.00000	0.00000	0.00000	0.00000
9	0.06667	0.66667	0.13333	0.13333	0.00000	0.00000	0.00000	0.00000	0.00000	0.00000	0.00000	0.00000	0.00000	0.00000	0.00000	0.00000	0.00000	0.00000	0.00000	0.00000	0.00000	0.00000	0.00000	0.00000	0.00000	0.00000	0.00000	0.00000	0.00000
10	0.02778	0.63889	0.11111	0.19444	0.00000	0.02778	0.00000	0.00000	0.00000	0.00000	0.00000	0.00000	0.00000	0.00000	0.00000	0.00000	0.00000	0.00000	0.00000	0.00000	0.00000	0.00000	0.00000	0.00000	0.00000	0.00000	0.00000	0.00000	0.00000
11	0.01176	0.55294	0.16471	0.22353	0.02353	0.02353	0.00000	0.00000	0.00000	0.00000	0.00000	0.00000	0.00000	0.00000	0.00000	0.00000	0.00000	0.00000	0.00000	0.00000	0.00000	0.00000	0.00000	0.00000	0.00000	0.00000	0.00000	0.00000	0.00000
12	0.00442	0.46460	0.15929	0.24779	0.04425	0.03982	0.00000	0.00442	0.00000	0.00000	0.00000	0.00000	0.00000	0.00000	0.00000	0.00000	0.00000	0.00000	0.00000	0.00000	0.00000	0.00000	0.00000	0.00000	0.00000	0.00000	0.00000	0.00000	0.03540
13	0.00173	0.38754	0.17301	0.25433	0.06228	0.05536	0.00692	0.00346	0.00000	0.00000	0.00000	0.00000	0.00000	0.00000	0.00000	0.00000	0.00000	0.00000	0.00000	0.00000	0.00000	0.00000	0.00000	0.00000	0.00000	0.00000	0.00000	0.00000	0.05536
14	0.00064	0.31804	0.17272	0.25685	0.07138	0.06310	0.01020	0.00701	0.00000	0.00064	0.00000	0.00000	0.00000	0.00000	0.00000	0.00000	0.00000	0.00000	0.00000	0.00000	0.00000	0.00000	0.00000	0.00000	0.00000	0.00000	0.00000	0.00000	0.09943
15	0.00024	0.26137	0.16403	0.25477	0.07730	0.06811	0.01697	0.01061	0.00094	0.00047	0.00000	0.00000	0.00000	0.00000	0.00000	0.00000	0.00000	0.00000	0.00000	0.00000	0.00000	0.00000	0.00000	0.00000	0.00000	0.00000	0.00000	0.00000	0.14518
16	0.00008	0.21042	0.15648	0.24755	0.08018	0.07335	0.02043	0.01367	0.00219	0.00110	0.00000	0.00008	0.00000	0.00000	0.00000	0.00000	0.00000	0.00000	0.00000	0.00000	0.00000	0.00000	0.00000	0.00000	0.00000	0.00000	0.00000	0.00000	0.19446
17	0.00003	0.16992	0.14515	0.23683	0.08265	0.07461	0.02308	0.01577	0.00369	0.00187	0.00018	0.00006	0.00000	0.00000	0.00000	0.00000	0.00000	0.00000	0.00000	0.00000	0.00000	0.00000	0.00000	0.00000	0.00000	0.00000	0.00000	0.00000	0.24616
18	0.00001	0.13609	0.13274	0.22465	0.08262	0.07493	0.02479	0.01767	0.00467	0.00261	0.00036	0.00016	0.00000	0.00001	0.00000	0.00000	0.00000	0.00000	0.00000	0.00000	0.00000	0.00000	0.00000	0.00000	0.00000	0.00000	0.00000	0.00000	0.29867
19	0.00000	0.10889	0.11986	0.21063	0.08154	0.07410	0.02593	0.01865	0.00557	0.00327	0.00068	0.00030	0.00002	0.00001	0.00000	0.00000	0.00000	0.00000	0.00000	0.00000	0.00000	0.00000	0.00000	0.00000	0.00000	0.00000	0.00000	0.00000	0.35054
20	0.00000	0.08666	0.10727	0.19592	0.07930	0.07275	0.02641	0.01934	0.00618	0.00384	0.00092	0.00045	0.00006	0.00002	0.00000	0.00000	0.00000	0.00000	0.00000	0.00000	0.00000	0.00000	0.00000	0.00000	0.00000	0.00000	0.00000	0.00000	0.40089
}\testdata

\begin{figure}
    \centering
\begin{tikzpicture}
\begin{axis}[
            ybar stacked,   
            ymin=0,         
            ymax = 1,
            xtick=data,     
            xticklabels from table={\testdata}{Label}, width=12.6cm, height = 7cm,bar width=0.43cm, enlarge x limits=0.05
]

\def\plotcommand#1{
    \addplot [fill=black!#1!blue!80!white] table [y=\s, meta=Label,x expr=\coordindex] {\testdata};
}

\addplot [fill=black!0!blue!80!white] table [y=1, meta=Label,x expr=\coordindex] {\testdata};

\foreach \s in {2,...,14}
{
\pgfmathparse{ln(\s+1)*100/3}

\expandafter\plotcommand\expandafter{\pgfmathresult}
}

\addplot [fill=orange!80!white!80!red] table [y=29, meta=Label,x expr=\coordindex] {\testdata};
\end{axis}
\node[align=left] (2) at (12.5,0.3) {have $\HM(2)$\\but not $\HM(1)$};
\node[align=left] (4) at (12.5,1.7) {have $\HM(4)$\\but not $\HM(3)$};
\node[align=left] (6) at (12.5,2.8) {have $\HM(6)$\\but not $\HM(5)$};
\node[align=left] (NL) at (12,4.7) {likely \\ NL-hard};
\draw (10.5,0.3) -- (2);
\draw (10.5,1.7) -- (4);
\draw (10.5,2.8) -- (6);
\draw (10.5,4.7) -- (NL);
\end{tikzpicture}
    \caption{Distribution of core trees in L.}
    \label{fig:HMtrees}
\end{figure}

\subsubsection{The Smallest Tree not Solved by Datalog} 
It turns out that every tree with at most 20 vertices which is not NP-hard can be solved by Datalog, thus confirming Conjecture~\ref{conj:bulin}. 
In fact, up to 20 vertices all trees that have 
Kearnes-Markovi\'{c}-McKenzie polymorphisms
either have a majority polymorphism or totally symmetric polymorphisms of all arities. The picture is however more complex for larger trees: there exists a tree which can be solved by Datalog but does not have a near-unanimity polymorphism (of any arity) and does not have totally symmetric polymorphisms of all arities, see \cite[Proposition 5.5]{BartoB13} for an example and \cite{Bulin18} for its solvability in Datalog.

\subsubsection{The Smallest Tree not Solved by Arc Consistency} 
\label{sect:NAC}
The smallest tree $\T$ that has no binary symmetric polymorphism has 19 vertices and is displayed in Figure~\ref{fig:withWNU3withoutWNU2}. It has a 3-wnu polymorphism, and even a majority polymorphism which satisfies $f(x,x,y)=f(x,y,x)=f(y,x,x)=x$ for all $x,y \in T$. 
Note that $\Csp(\T)$ cannot be solved by the arc-consistency procedure since in this case $\T$ must have a binary symmetric polymorphism~\cite{FederVardi,DalmauPearson}. 
All other trees with at most 19 vertices satisfy $\TS(n)$ for all $n$. For a tree $\T$ the vertices of the indicator digraph for $\TS(2 |E(\T)|)$ correspond to the nonempty subsets of $T$. Hence the indicator structure of a tree $\T$ with 19 vertices has $2^{19}-1=524287$ vertices. Using level-wise satisfiability (see Section~\ref{sect:AC-Pol}) the number of vertices of the indicator structure is reduced to something between 19 and 513, depending on the number of vertices on each level. 

\begin{figure}
\centering\begin{tikzpicture}[scale=0.4]
\node[bullet] (4) at (0,0) {};
\node[bullet] (18) at (1,0) {};
\node[bullet] (3) at (0,1) {};
\node[bullet] (5) at (4,1) {};
\node[bullet] (11) at (3,1) {};
\node[bullet] (14) at (2,1) {};
\node[bullet] (17) at (1,1) {};
\node[bullet] (2) at (0,2) {};
\node[bullet] (15) at (1,2) {};
\node[bullet] (0) at (2,2) {};
\node[bullet] (6) at (4,2) {};
\node[bullet] (10) at (3,2) {};
\node[bullet] (1) at (0,3) {};
\node[bullet] (7) at (4,3) {};
\node[bullet] (9) at (2,3) {};
\node[bullet] (12) at (3,3) {};
\node[bullet] (16) at (1,3) {};
\node[bullet] (8) at (4,4) {};
\node[bullet] (13) at (3,4) {};
\path[->,>=stealth']
(0) edge (1)
(0) edge (9)
(2) edge (1)
(3) edge (2)
(4) edge (3)
(5) edge (0)
(5) edge (6)
(6) edge (7)
(7) edge (8)
(10) edge (9)
(10) edge (12)
(11) edge (10)
(12) edge (13)
(14) edge (0)
(14) edge (15)
(15) edge (16)
(18) edge (17)
(17) edge (15)
;

\node at (2,-1) {Tree C};
\end{tikzpicture}
\caption{The smallest tree that cannot be solved by Arc Consistency (it has 19 vertices and a majority, but no binary symmetric polymorphism).}
\label{fig:withWNU3withoutWNU2}
\end{figure}

\subsection{Open Trees}
\label{sect:opentrees}
In this section we present trees that are interesting test cases, 
in particular regarding the conjectured classification of digraphs in NL (Conjecture~\ref{conj:NL}).

\subsubsection{A Tree not Known to be in NL}
We found a tree with polymorphisms that form a Kearnes-Kiss chain of length five and a Hobby-McKenzie chain of length 2, but has no Majority, and no (level-wise) J\'onsson chain of length 1000 (see Figure~\ref{fig:KKtree}). This tree is neither known to be P-hard or Mod$_p$L-hard, nor is it known to be in NL.  
It is the smallest tree without a majority polymorphism. 
Note that the existence of a J\'onsson chain of \emph{some} length is decidable because for a given digraph there are only finitely many operations of arity three. Moreover, by the discussion from Section~\ref{sect:level-wise} we know that we may narrow down the set of operations that have to be considered; the resulting number of operations is $12^{36}$. Even if we could show that the tree has no J\'onsson chain
we would not know that the tree is not in NL. We believe that the tree is in NL, but new ideas are needed to prove that (e.g., ideas to prove Conjecture~\ref{conj:NL}).  We mention that it can pp-construct $\stcon$, so it is NL-hard.

\def\scale{0.4}
\def\hdist{1cm}
\def\vdist{1cm}
\begin{figure}
\centering\begin{tikzpicture}[scale=\scale]
\node[bullet,gray] (12) at (1.5,0) {};
\node[bullet] (1) at (0.5,1) {};
\node[bullet] (10) at (3,1) {};
\node[bullet] (11) at (1.5,1) {};
\node[bullet] (2) at (0,2) {};
\node[bullet] (0) at (1,2) {};
\node[bullet] (7) at (3,2) {};
\node[bullet] (13) at (2,2) {};
\node[bullet] (3) at (0,3) {};
\node[bullet] (6) at (1,3) {};
\node[bullet] (8) at (3,3) {};
\node[bullet] (14) at (2,3) {};
\node[bullet] (4) at (0,4) {};
\node[bullet] (9) at (3,4) {};
\node[bullet] (15) at (2,4) {};
\node[bullet,gray] (5) at (0,5) {};
\path[->,>=stealth']
(1) edge (2)
(1) edge (0)
(2) edge (3)
(3) edge (4)
(4) edge[gray] (5)
(0) edge (6)
(7) edge (6)
(7) edge (8)
(10) edge (7)
(8) edge (9)
(11) edge (0)
(11) edge (13)
(12) edge[gray] (11)
(13) edge (14)
(14) edge (15)
;
\node at (1.5,-1) {Tree D};
\end{tikzpicture}
\caption{The smallest tree without a majority polymorphism (unique up to edge reversal; it has 16 vertices). 
It satisfies KK(5) but does not satisfy J\'onsson(1000). Therefore, Conjecture~\ref{conj:NL} puts it in NL but we cannot prove this fact; it is an interesting open case. 
}
\label{fig:KKtree}
\end{figure}

\subsubsection{Trees that might be P-hard}
There are 28 trees with 18 vertices that satisfy neither $\HMcK(1000)$ nor $\KK(1000)$, not even level-wise (see Figure~\ref{fig:noKKtrees}). 
They satisfy $\TS(n)$ for all $n$, so they are in P and cannot pp-construct $\TLinP$ for every $p$. Hence, this is in accordance with Conjecture~\ref{conj:bulin}. All other trees with up to 18 vertices satisfy $\KK(5)$ and are in NL assuming Conjecture~\ref{conj:NL}. 
Hence, if this conjecture is true, and if $\operatorname{NL} \neq \operatorname{P}$, and if indeed
these 28 trees do not have $\HMcK(n)$ for any $n$, then they are the smallest trees that are P-hard.

\def\scale{0.4}
\def\hdist{1cm}
\def\vdist{1cm}
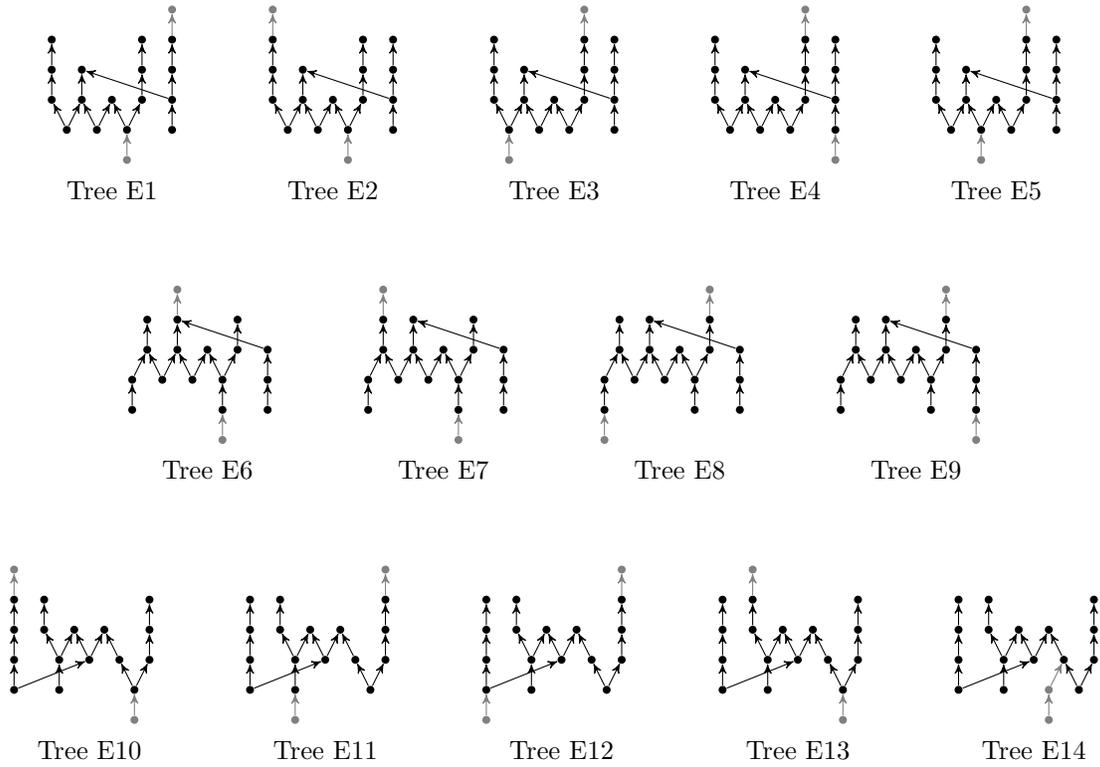
\begin{figure}[h]
\centering
\begin{tikzpicture}[scale=\scale]
\node[bullet,gray] (3) at (2.5,0) {};
\node[bullet] (2) at (2.5,1) {};
\node[bullet] (14) at (4,1) {};
\node[bullet] (0) at (1.5,1) {};
\node[bullet] (8) at (0.5,1) {};
\node[bullet] (1) at (2,2) {};
\node[bullet] (13) at (4,2) {};
\node[bullet] (7) at (1,2) {};
\node[bullet] (9) at (0,2) {};
\node[bullet] (4) at (3,2) {};
\node[bullet] (5) at (3,3) {};
\node[bullet] (10) at (0,3) {};
\node[bullet] (15) at (4,3) {};
\node[bullet] (12) at (1,3) {};
\node[bullet] (6) at (3,4) {};
\node[bullet] (11) at (0,4) {};
\node[bullet] (16) at (4,4) {};
\node[bullet,gray] (17) at (4,5) {};
\path[->,>=stealth']
(2) edge (1)
(2) edge (4)
(5) edge (6)
(14) edge (13)
(13) edge (15)
(13) edge (12)
(10) edge (11)
(0) edge (7)
(0) edge (1)
(7) edge (12)
(8) edge (9)
(8) edge (7)
(9) edge (10)
(4) edge (5)
(3) edge[gray] (2)
(15) edge (16)
(16) edge[gray] (17)
;
\node at (2,-1) {Tree E1};
\end{tikzpicture}
\hspace{\hdist}
\begin{tikzpicture}[scale=\scale]
\node[bullet,gray] (3) at (2.5,0) {};
\node[bullet] (13) at (0.5,1) {};
\node[bullet] (2) at (2.5,1) {};
\node[bullet] (10) at (4,1) {};
\node[bullet] (0) at (1.5,1) {};
\node[bullet] (4) at (3,2) {};
\node[bullet] (7) at (1,2) {};
\node[bullet] (9) at (4,2) {};
\node[bullet] (14) at (0,2) {};
\node[bullet] (1) at (2,2) {};
\node[bullet] (5) at (3,3) {};
\node[bullet] (15) at (0,3) {};
\node[bullet] (8) at (1,3) {};
\node[bullet] (11) at (4,3) {};
\node[bullet] (16) at (0,4) {};
\node[bullet] (6) at (3,4) {};
\node[bullet] (12) at (4,4) {};
\node[bullet,gray] (17) at (0,5) {};
\path[->,>=stealth']
(4) edge (5)
(5) edge (6)
(13) edge (7)
(13) edge (14)
(7) edge (8)
(3) edge[gray] (2)
(2) edge (4)
(2) edge (1)
(15) edge (16)
(16) edge[gray] (17)
(10) edge (9)
(9) edge (8)
(9) edge (11)
(14) edge (15)
(0) edge (7)
(0) edge (1)
(11) edge (12)
;
\node at (2,-1) {Tree E2};
\end{tikzpicture}
\hspace{\hdist}
\begin{tikzpicture}[scale=\scale]
\node[bullet,gray] (9) at (0.5,0) {};
\node[bullet] (1) at (2.5,1) {};
\node[bullet] (8) at (0.5,1) {};
\node[bullet] (6) at (1.5,1) {};
\node[bullet] (15) at (4,1) {};
\node[bullet] (2) at (3,2) {};
\node[bullet] (10) at (0,2) {};
\node[bullet] (0) at (2,2) {};
\node[bullet] (14) at (4,2) {};
\node[bullet] (7) at (1,2) {};
\node[bullet] (3) at (3,3) {};
\node[bullet] (13) at (1,3) {};
\node[bullet] (11) at (0,3) {};
\node[bullet] (16) at (4,3) {};
\node[bullet] (4) at (3,4) {};
\node[bullet] (17) at (4,4) {};
\node[bullet] (12) at (0,4) {};
\node[bullet,gray] (5) at (3,5) {};
\path[->,>=stealth']
(3) edge (4)
(4) edge[gray] (5)
(1) edge (2)
(1) edge (0)
(2) edge (3)
(9) edge[gray] (8)
(8) edge (10)
(8) edge (7)
(10) edge (11)
(6) edge (0)
(6) edge (7)
(14) edge (13)
(14) edge (16)
(11) edge (12)
(16) edge (17)
(7) edge (13)
(15) edge (14)
;
\node at (2,-1) {Tree E3};
\end{tikzpicture}
\hspace{\hdist}
\begin{tikzpicture}[scale=\scale]
\node[bullet,gray] (15) at (4,0) {};
\node[bullet] (14) at (4,1) {};
\node[bullet] (6) at (1.5,1) {};
\node[bullet] (8) at (0.5,1) {};
\node[bullet] (1) at (2.5,1) {};
\node[bullet] (7) at (1,2) {};
\node[bullet] (0) at (2,2) {};
\node[bullet] (13) at (4,2) {};
\node[bullet] (9) at (0,2) {};
\node[bullet] (2) at (3,2) {};
\node[bullet] (10) at (0,3) {};
\node[bullet] (12) at (1,3) {};
\node[bullet] (3) at (3,3) {};
\node[bullet] (16) at (4,3) {};
\node[bullet] (11) at (0,4) {};
\node[bullet] (4) at (3,4) {};
\node[bullet] (17) at (4,4) {};
\node[bullet,gray] (5) at (3,5) {};
\path[->,>=stealth']
(15) edge[gray] (14)
(14) edge (13)
(10) edge (11)
(7) edge (12)
(6) edge (0)
(6) edge (7)
(13) edge (12)
(13) edge (16)
(8) edge (9)
(8) edge (7)
(9) edge (10)
(3) edge (4)
(4) edge[gray] (5)
(1) edge (2)
(1) edge (0)
(2) edge (3)
(16) edge (17)
;
\node at (2,-1) {Tree E4};
\end{tikzpicture}
\hspace{\hdist}
\begin{tikzpicture}[scale=\scale]
\node[bullet,gray] (7) at (1.5,0) {};
\node[bullet] (9) at (0.5,1) {};
\node[bullet] (1) at (2.5,1) {};
\node[bullet] (6) at (1.5,1) {};
\node[bullet] (15) at (4,1) {};
\node[bullet] (10) at (0,2) {};
\node[bullet] (0) at (2,2) {};
\node[bullet] (8) at (1,2) {};
\node[bullet] (14) at (4,2) {};
\node[bullet] (2) at (3,2) {};
\node[bullet] (16) at (4,3) {};
\node[bullet] (13) at (1,3) {};
\node[bullet] (11) at (0,3) {};
\node[bullet] (3) at (3,3) {};
\node[bullet] (17) at (4,4) {};
\node[bullet] (12) at (0,4) {};
\node[bullet] (4) at (3,4) {};
\node[bullet,gray] (5) at (3,5) {};
\path[->,>=stealth']
(16) edge (17)
(9) edge (10)
(9) edge (8)
(10) edge (11)
(1) edge (0)
(1) edge (2)
(8) edge (13)
(6) edge (0)
(6) edge (8)
(11) edge (12)
(7) edge[gray] (6)
(4) edge[gray] (5)
(14) edge (16)
(14) edge (13)
(15) edge (14)
(3) edge (4)
(2) edge (3)
;
\node at (2,-1) {Tree E5};
\end{tikzpicture}

\vspace{\vdist}

\begin{tikzpicture}[scale=\scale]
\node[bullet,gray] (4) at (2.5,0) {};
\node[bullet] (3) at (2.5,1) {};
\node[bullet] (17) at (-0.5,1) {};
\node[bullet] (12) at (4,1) {};
\node[bullet] (11) at (4,2) {};
\node[bullet] (0) at (1.5,2) {};
\node[bullet] (13) at (0.5,2) {};
\node[bullet] (16) at (-0.5,2) {};
\node[bullet] (2) at (2.5,2) {};
\node[bullet] (10) at (4,3) {};
\node[bullet] (7) at (1,3) {};
\node[bullet] (14) at (0,3) {};
\node[bullet] (1) at (2,3) {};
\node[bullet] (5) at (3,3) {};
\node[bullet] (8) at (1,4) {};
\node[bullet] (15) at (0,4) {};
\node[bullet] (6) at (3,4) {};
\node[bullet,gray] (9) at (1,5) {};
\path[->,>=stealth']
(4) edge[gray] (3)
(3) edge (2)
(11) edge (10)
(10) edge (8)
(0) edge (7)
(0) edge (1)
(7) edge (8)
(8) edge[gray] (9)
(14) edge (15)
(13) edge (14)
(13) edge (7)
(16) edge (14)
(17) edge (16)
(2) edge (1)
(2) edge (5)
(5) edge (6)
(12) edge (11)
;
\node at (2,-1) {Tree E6};
\end{tikzpicture}
\hspace{\hdist}
\begin{tikzpicture}[scale=\scale]
\node[bullet,gray] (4) at (2.5,0) {};
\node[bullet] (15) at (-0.5,1) {};
\node[bullet] (11) at (4,1) {};
\node[bullet] (3) at (2.5,1) {};
\node[bullet] (0) at (1.5,2) {};
\node[bullet] (14) at (-0.5,2) {};
\node[bullet] (2) at (2.5,2) {};
\node[bullet] (10) at (4,2) {};
\node[bullet] (12) at (0.5,2) {};
\node[bullet] (5) at (3,3) {};
\node[bullet] (1) at (2,3) {};
\node[bullet] (13) at (0,3) {};
\node[bullet] (7) at (1,3) {};
\node[bullet] (9) at (4,3) {};
\node[bullet] (6) at (3,4) {};
\node[bullet] (16) at (0,4) {};
\node[bullet] (8) at (1,4) {};
\node[bullet,gray] (17) at (0,5) {};
\path[->,>=stealth']
(5) edge (6)
(0) edge (1)
(0) edge (7)
(15) edge (14)
(14) edge (13)
(2) edge (5)
(2) edge (1)
(13) edge (16)
(16) edge[gray] (17)
(11) edge (10)
(10) edge (9)
(4) edge[gray] (3)
(3) edge (2)
(12) edge (7)
(12) edge (13)
(7) edge (8)
(9) edge (8)
;
\node at (2,-1) {Tree E7};
\end{tikzpicture}
\hspace{\hdist}
\begin{tikzpicture}[scale=\scale]
\node[bullet,gray] (17) at (-0.5,0) {};
\node[bullet] (11) at (4,1) {};
\node[bullet] (3) at (2.5,1) {};
\node[bullet] (16) at (-0.5,1) {};
\node[bullet] (10) at (4,2) {};
\node[bullet] (12) at (0.5,2) {};
\node[bullet] (2) at (2.5,2) {};
\node[bullet] (0) at (1.5,2) {};
\node[bullet] (15) at (-0.5,2) {};
\node[bullet] (7) at (1,3) {};
\node[bullet] (1) at (2,3) {};
\node[bullet] (4) at (3,3) {};
\node[bullet] (9) at (4,3) {};
\node[bullet] (13) at (0,3) {};
\node[bullet] (8) at (1,4) {};
\node[bullet] (5) at (3,4) {};
\node[bullet] (14) at (0,4) {};
\node[bullet,gray] (6) at (3,5) {};
\path[->,>=stealth']
(11) edge (10)
(10) edge (9)
(7) edge (8)
(12) edge (7)
(12) edge (13)
(3) edge (2)
(2) edge (1)
(2) edge (4)
(0) edge (1)
(0) edge (7)
(4) edge (5)
(5) edge[gray] (6)
(17) edge[gray] (16)
(16) edge (15)
(9) edge (8)
(15) edge (13)
(13) edge (14)
;
\node at (2,-1) {Tree E8};
\end{tikzpicture}
\hspace{\hdist}
\begin{tikzpicture}[scale=\scale]
\node[bullet,gray] (17) at (4,0) {};
\node[bullet] (16) at (4,1) {};
\node[bullet] (12) at (-0.5,1) {};
\node[bullet] (3) at (2.5,1) {};
\node[bullet] (0) at (1.5,2) {};
\node[bullet] (11) at (-0.5,2) {};
\node[bullet] (8) at (0.5,2) {};
\node[bullet] (2) at (2.5,2) {};
\node[bullet] (15) at (4,2) {};
\node[bullet] (7) at (1,3) {};
\node[bullet] (9) at (0,3) {};
\node[bullet] (1) at (2,3) {};
\node[bullet] (4) at (3,3) {};
\node[bullet] (14) at (4,3) {};
\node[bullet] (13) at (1,4) {};
\node[bullet] (5) at (3,4) {};
\node[bullet] (10) at (0,4) {};
\node[bullet,gray] (6) at (3,5) {};
\path[->,>=stealth']
(0) edge (7)
(0) edge (1)
(7) edge (13)
(11) edge (9)
(9) edge (10)
(8) edge (7)
(8) edge (9)
(2) edge (1)
(2) edge (4)
(4) edge (5)
(15) edge (14)
(14) edge (13)
(17) edge[gray] (16)
(16) edge (15)
(12) edge (11)
(5) edge[gray] (6)
(3) edge (2)
;
\node at (2,-1) {Tree E9};
\end{tikzpicture}

\vspace{\vdist}

\begin{tikzpicture}[scale=\scale]
\node[bullet,gray] (3) at (3.5,0) {};
\node[bullet] (13) at (-0.5,1) {};
\node[bullet] (2) at (3.5,1) {};
\node[bullet] (10) at (1,1) {};
\node[bullet] (9) at (1,2) {};
\node[bullet] (1) at (3,2) {};
\node[bullet] (7) at (2,2) {};
\node[bullet] (14) at (-0.5,2) {};
\node[bullet] (4) at (4,2) {};
\node[bullet] (8) at (1.5,3) {};
\node[bullet] (0) at (2.5,3) {};
\node[bullet] (5) at (4,3) {};
\node[bullet] (11) at (0.5,3) {};
\node[bullet] (15) at (-0.5,3) {};
\node[bullet] (12) at (0.5,4) {};
\node[bullet] (16) at (-0.5,4) {};
\node[bullet] (6) at (4,4) {};
\node[bullet,gray] (17) at (-0.5,5) {};
\path[->,>=stealth']
(9) edge (8)
(9) edge (11)
(1) edge (0)
(7) edge (8)
(7) edge (0)
(13) edge (14)
(13) edge (7)
(14) edge (15)
(4) edge (5)
(5) edge (6)
(3) edge[gray] (2)
(2) edge (1)
(2) edge (4)
(10) edge (9)
(11) edge (12)
(15) edge (16)
(16) edge[gray] (17)
;
\node at (2,-1) {Tree E10};
\end{tikzpicture}
\hspace{\hdist}
\begin{tikzpicture}[scale=\scale]
\node[bullet,gray] (15) at (1,0) {};
\node[bullet] (8) at (-0.5,1) {};
\node[bullet] (1) at (3.5,1) {};
\node[bullet] (14) at (1,1) {};
\node[bullet] (13) at (1,2) {};
\node[bullet] (0) at (3,2) {};
\node[bullet] (9) at (-0.5,2) {};
\node[bullet] (2) at (4,2) {};
\node[bullet] (7) at (2,2) {};
\node[bullet] (16) at (0.5,3) {};
\node[bullet] (6) at (2.5,3) {};
\node[bullet] (10) at (-0.5,3) {};
\node[bullet] (12) at (1.5,3) {};
\node[bullet] (3) at (4,3) {};
\node[bullet] (17) at (0.5,4) {};
\node[bullet] (11) at (-0.5,4) {};
\node[bullet] (4) at (4,4) {};
\node[bullet,gray] (5) at (4,5) {};
\path[->,>=stealth']
(13) edge (16)
(13) edge (12)
(16) edge (17)
(0) edge (6)
(8) edge (9)
(8) edge (7)
(9) edge (10)
(10) edge (11)
(1) edge (2)
(1) edge (0)
(2) edge (3)
(7) edge (6)
(7) edge (12)
(14) edge (13)
(4) edge[gray] (5)
(3) edge (4)
(15) edge[gray] (14)
;
\node at (2,-1) {Tree E11};
\end{tikzpicture}
\hspace{\hdist}
\begin{tikzpicture}[scale=\scale]
\node[bullet,gray] (9) at (-0.5,0) {};
\node[bullet] (1) at (3.5,1) {};
\node[bullet] (8) at (-0.5,1) {};
\node[bullet] (15) at (1,1) {};
\node[bullet] (2) at (4,2) {};
\node[bullet] (7) at (2,2) {};
\node[bullet] (0) at (3,2) {};
\node[bullet] (14) at (1,2) {};
\node[bullet] (10) at (-0.5,2) {};
\node[bullet] (3) at (4,3) {};
\node[bullet] (6) at (2.5,3) {};
\node[bullet] (11) at (-0.5,3) {};
\node[bullet] (16) at (0.5,3) {};
\node[bullet] (13) at (1.5,3) {};
\node[bullet] (4) at (4,4) {};
\node[bullet] (12) at (-0.5,4) {};
\node[bullet] (17) at (0.5,4) {};
\node[bullet,gray] (5) at (4,5) {};
\path[->,>=stealth']
(2) edge (3)
(3) edge (4)
(4) edge[gray] (5)
(7) edge (6)
(7) edge (13)
(11) edge (12)
(1) edge (0)
(1) edge (2)
(0) edge (6)
(8) edge (7)
(8) edge (10)
(14) edge (16)
(14) edge (13)
(16) edge (17)
(10) edge (11)
(15) edge (14)
(9) edge[gray] (8)
;
\node at (2,-1) {Tree E12};
\end{tikzpicture}
\hspace{\hdist}
\begin{tikzpicture}[scale=\scale]
\node[bullet,gray] (3) at (3.5,0) {};
\node[bullet] (2) at (3.5,1) {};
\node[bullet] (14) at (1,1) {};
\node[bullet] (8) at (-0.5,1) {};
\node[bullet] (4) at (4,2) {};
\node[bullet] (1) at (3,2) {};
\node[bullet] (13) at (1,2) {};
\node[bullet] (9) at (-0.5,2) {};
\node[bullet] (7) at (2,2) {};
\node[bullet] (5) at (4,3) {};
\node[bullet] (10) at (-0.5,3) {};
\node[bullet] (15) at (0.5,3) {};
\node[bullet] (0) at (2.5,3) {};
\node[bullet] (12) at (1.5,3) {};
\node[bullet] (6) at (4,4) {};
\node[bullet] (11) at (-0.5,4) {};
\node[bullet] (16) at (0.5,4) {};
\node[bullet,gray] (17) at (0.5,5) {};
\path[->,>=stealth']
(5) edge (6)
(4) edge (5)
(2) edge (1)
(2) edge (4)
(1) edge (0)
(14) edge (13)
(13) edge (12)
(13) edge (15)
(10) edge (11)
(15) edge (16)
(16) edge[gray] (17)
(9) edge (10)
(8) edge (7)
(8) edge (9)
(7) edge (12)
(7) edge (0)
(3) edge[gray] (2)
;
\node at (2,-1) {Tree E13};
\end{tikzpicture}
\hspace{\hdist}
\begin{tikzpicture}[scale=\scale]
\node[bullet] (12) at (-0.5,4) {};
\node[bullet] (17) at (0.5,4) {};
\node[bullet] (7) at (4,4) {};
\node[bullet] (6) at (4,3) {};
\node[bullet] (11) at (-0.5,3) {};
\node[bullet] (16) at (0.5,3) {};
\node[bullet] (0) at (2.5,3) {};
\node[bullet] (13) at (1.5,3) {};
\node[bullet] (5) at (4,2) {};
\node[bullet] (14) at (1,2) {};
\node[bullet] (8) at (2,2) {};
\node[bullet] (1) at (3,2) {};
\node[bullet] (10) at (-0.5,2) {};
\node[bullet,gray] (2) at (2.5,1) {};
\node[bullet] (15) at (1,1) {};
\node[bullet] (4) at (3.5,1) {};
\node[bullet] (9) at (-0.5,1) {};
\node[bullet,gray] (3) at (2.5,0) {};
\path[->,>=stealth']
(5)  edge (6)
(4)  edge (5)
(11) edge (12)
(10) edge (11)
(3)  edge[gray] (2)
(15) edge (14)
(16) edge (17)
(14) edge (16)
(8)  edge (0)
(1)  edge (0)
(9)  edge (8)
(4)  edge (1)
(2)  edge[gray] (1)
(9)  edge (10)
(14) edge (13)
(6)  edge (7)
(8)  edge (13)
;
\node at (2,-1) {Tree E14};
\end{tikzpicture}
\caption{Smallest trees without $\HMcK(1000)$ (they have 18 vertices). These trees are tractable and candidates for being P-complete.
}
\label{fig:noKKtrees}
\end{figure}

\subsection{Majority Polymorphisms}
\label{sect:other} 
Majority polymorphisms play a central role in the early theory of the constraint satisfaction problem~\cite{FederVardi,FederCycles,JeavonsClosure},
in graph theory~\cite{Kazda,HellRafiey-list-homomorphism-digraphs},
and in the algebraic theory of CSPs~\cite{Bulatov-Conservative-Revisited,BulatovFVConjecture}. 
We have therefore also computed a
smallest tree without a majority polymorphism (see Figure~\ref{fig:KKtree}). Interestingly, when solving the indicator problem for the existence of a majority  polymorphism of $\HH$ for graphs with at most 15 vertices (which all have a majority polymorphism), 
no backtracking was needed: pruning with the arc-consistency procedure did suffice to avoid all dead-ends in the search. Theoretical results only guarantee this behavior for establishing $(2,3)$-consistency (since $\HH$ has a majority polymorphism). 
So one might ask: can every tree with a majority polymorphism be solved by arc consistency? 
This is not the case; see Lemmata 4.1 and 4.2 in~\cite{SpecialTriads}. 
In our experiments we found the smallest such tree: 
Figure~\ref{fig:withWNU3withoutWNU2} shows a tree with a majority polymorphism which does
not even have a binary symmetric polymorphism, and hence in particular cannot be solved by arc consistency.

\section{Open Problems and Future Work} \label{sect:discussion}
The following conjecture is implied by 
Conjecture~\ref{conj:bulin}, but might be easier to answer.

\begin{conjecture}
A tree has Kearnes-Markovi\'{c}-McKenzie polymorphisms if and only if it has a 3-wnu polymorphism.
\end{conjecture}

\begin{question}
Is it true that the probability that a tree drawn uniformly at random from the set of all trees with vertex set $\{1,\dots,n\}$ is NP-hard tends to 1 as $n$ tends to infinity? The answer is yes if we ask the question for random labelled digraphs instead of random labelled trees~\cite{LuczakNesetril}. 
\end{question}

The table in Figure~\ref{fig:HMtrees} suggests that the following conjecture is true.

\begin{conjecture}
The fraction of core trees with $n$ vertices that are NL-hard goes to 1 as $n$ goes to infinity.
\end{conjecture}

 \begin{question}
 Determine the smallest trees that are P-hard (assuming that $\operatorname{NL}\neq\operatorname{P}$). We know from Section~\ref{sect:other} that they must have at least 16 vertices, since all smaller trees have a majority and thus are in NL. 
\end{question}

 \begin{question}
Is our algorithm from Section~\ref{sect:gen} to generate unlabeled core trees a polynomial-delay enumeration algorithm (in the sense of~\cite{JohnYannaPapaGen})?
\end{question}
 
\begin{question}
Characterize linear conditions that can be tested level-wise (in the sense of Section~\ref{sect:level-wise}) for balanced digraphs, and more specifically, for trees.
\end{question}

 It would be interesting to perform experiments similar to the experiments presented here for trees that are equipped with a singleton unary relation $\{a\}$ for each vertex $a$ of the tree; in this case, if $\T_c$ is the resulting expanded tree structure, $\Csp(\T_c)$ models the so-called \emph{$\T$-precoloring extension problem}. This setting is particularly nice from the algebraic perspective because 
 then all the polymorphisms of $\T_c$ 
 are idempotent. Note, however, that all these structure $\T_c$ are cores, so there are far more structures to consider, and hardness will dominate more rapidly.

Taking this one step further, it would also be interesting to study the so-called \emph{list homomorphism problem} for trees $\T$ from an experimental perspective. Here, the input contains besides the graph $\G$ a list of vertices from $\HH$ and we are looking for a homomorphism from $\G$ to $\HH$ that maps each vertex to an element from its list. This can be seen as a special case of a CSP for a relational structure, which contains besides the edge relation also a unary relation for each subset of the vertices of $\HH$. On the algebraic side, we are therefore interested in polymorphisms that preserve all subsets of $\HH$; such polymorphisms (and consequently the respective CSPs) are also called \emph{conservative}. The algorithms and complexities for conservative CSPs are better understood than the general case~\cite{Conservative, Barto-Conservative,Bulatov-Conservative-Revisited,Kazda-binary-conservative}, which will help to determine the complexity of the list homomorphism for trees. On the other hand, as in the case of the precoloring extension problem we have much larger numbers of trees to consider since all the structures that we study are already cores.

\section{Declarations}
The authors are grateful to the Center for Information Services and High Performance Computing [Zentrum f\"ur Informationsdienste und Hochleistungsrechnen (ZIH)] at TU Dresden for providing its facilities for high throughput calculations.
Manuel Bodirsky has received funding from the European Research Council (Grant Agreement no. 681988, CSP-Infinity).
Jakub Bul{\' i}n was supported by the M{\v S}MT {\v C}R INTER-EXCELLENCE project LTAUSA19070 and the Charles University project UNCE/SCI/004.
Florian Starke is supported by DFG Graduiertenkolleg 1763 (QuantLA).
The authors have no competing interests to declare that are relevant to the content of this article.

\bibliographystyle{abbrv}
\bibliography{global}

\end{document}